\documentstyle[12pt,leqno,twoside]{article}
\input{amssym.def}

\title{\LARGE \rm Special almost Hermitian geometry}
\date{}

\newcounter{sec}
\newtheorem{proposition}{Proposition}[sec]
\newtheorem{lemma}[proposition]{Lemma}

\newtheorem{corollary}[proposition]{Corollary}
\newtheorem{theorem}[proposition]{Theorem}
\newtheorem{observation}[proposition]{\it Remark}
\newcommand{\Dt}{{\it Proof}.- }

\newcommand{\com}{\makebox[7pt]{\raisebox{1.5pt}{\tiny$\circ$}}}

\newcommand{\cqd}{\hfill $\Box$\par}
\newcommand{\seccion}[1]{\setcounter{equation}{0}\addtocounter{sec}{1}
   \section{\centerline{\protect\large\bf \arabic{sec}. #1
\hspace{3em} }} \indent}

\newcommand{\SP}{\mathop{\mbox{\sl Sp}}}

\newcommand{\SUn}{\mathop{\mbox{\sl SU}}}

\newcommand{\lcf}{\lbrack\!\lbrack}
\newcommand{\rcf}{\rbrack\!\rbrack}

\def\thebibliography#1{\section*{\large\bf
References}
        \list {\bf \arabic{enumi}.}{\settowidth\labelwidth{[#1]}
 \leftmargin\labelwidth \advance\leftmargin\labelsep
                 \usecounter{enumi}}
         \def\newblock{\hskip .11em plus .33em minus .07em}
         \sloppy\clubpenalty4000\widowpenalty4000
         \sfcode`\.=1000\relax}

\begin{document}

\baselineskip=5mm \pagestyle{myheadings}

\maketitle
\vspace{-1cm}

\noindent FRANCISCO MART{\'I}N CABRERA

\begin{footnotesize}
  \setlength{\parindent}{0pt} \setlength{\parskip}{10pt}

\textit{Departamento de Matem{\'a}tica Fundamental, Universidad de La
Laguna, 38200 La Laguna, Tenerife, Spain}. E-mail: {\tt
fmartin@ull.es} \vspace{7mm}

{\bf Abstract.} We  study the classification of special almost
hermitian manifolds  in Gray and Hervella's type classes. We prove
that the exterior derivatives of the K{\"a}hler form and the
complex volume form contain all the information about the
intrinsic torsion of the $\SUn(n)$-structure. Furthermore, we
apply the obtained results to almost hyperhermitian geometry.
Thus, we show that the exterior derivatives of the three
K{\"a}hler forms of an almost hyperhermitian manifold are
sufficient to determine the three covariant derivatives of such
forms, i.e., the three mentioned  exterior derivatives determine
the intrinsic torsion of the  ${\sl Sp}(n)$-structure.

\vspace{7mm}

  \vspace{4mm}


\end{footnotesize}

\baselineskip=5mm

\seccion{Introduction}{\indent} \setcounter{proposition}{0} In
1955, Berger \cite{Berger} gave the list of possible holonomy
groups of non-symmetric Riemannian $m$-manifolds whose holonomy
representation is irreducible. Such a list of groups  was
complemented with their corresponding holonomy representations,
i.e., it was also specified the action of each group on the
tangent space. Consequently, each group $G \subseteq {\sl SO}(m)$
in Berger's list gives rise to a geometric structure. Moreover,
the groups $G$ may be given as the stabilisers in ${\sl SO}(m)$ of
certain differential forms on ${\Bbb R}^m$. For $G={\sl G}_2$, it
is a three-form $\phi$ on ${\Bbb R}^7$; for $G={\sl Spin}(7)$, it
is a four-form $\varphi$ on ${\Bbb R}^{8}$; for $G={\sl
Sp}(n)Sp(1)$, it is a four-form $\Omega$ on ${\Bbb R}^{4n}$;  for
$G={\sl U}(n)$,  a K{\"a}hler form $\omega$ on ${\Bbb R}^{2n}$, etc.
Such forms are a key ingredient in the definition of the
corresponding $G$-structure on a Riemannian $m$-manifold $M$.
Furthermore, the intrinsic torsion of a $G$-structure, defined in
next section, can be identified with the Levi-Civita covariant
derivatives of the corresponding forms and is always contained in
${\cal W} = T^* M \otimes \frak g^{\perp}$, being ${\frak s}{\frak
o}(m) = {\frak g} \oplus {\frak g}^{\perp}$. The action of $G$
splits ${\cal W}$ into irreducible components, say ${\cal W} =
{\cal W}_1 \oplus \ldots \oplus {\cal W}_k$. Then, $G$-structures
on $M$ can be classified in at most $2^k$ classes.

This way of classifying $G$-structures was initiated by Gray and
Hervella in \cite{GrayHervella}, where they considered the case $G
= {\sl U}(n)$ (almost Hermitian structures), turning out ${\cal W}
= {\cal W}_1 \oplus {\cal W}_2 \oplus {\cal W}_3 \oplus {\cal
W}_4$, for $n
>2$, i.e., there are sixteen classes of almost Hermitian
manifolds. Later, diverse authors have studied the situation for
other $G$-structures: ${\sl G}_2$, ${\sl Spin}(7)$, ${\sl
Sp}(n){\sl Sp}(1)$, etc.

In the present paper we study the situation for $G={\sl SU}(n)$.
Thus, we consider Riemannian $2n$-manifolds equipped with a K{\"a}hler
form $\omega$ and a complex volume form $\Psi = \psi_+ + i
\psi_-$, called special almost Hermitian manifolds. The group
${\sl SU}(n)$ is the stabiliser in ${\sl SO}(2n)$ of $\omega$ and
$\Psi$. Therefore, the information about intrinsic torsion of an
${\sl SU}(n)$-structure is contained in $\nabla \omega$ and
$\nabla \Psi$, where $\nabla$ denotes the Levi-Civita connection.
For high dimensions, $2n \geq 8$, we find
$$
T^* M \otimes {\frak s}{\frak u}(n)^{\perp} = {\cal W}_1 \oplus
{\cal W}_2 \oplus {\cal W}_3 \oplus {\cal W}_4 \oplus {\cal W}_5,
$$
where the first four summands coincide with Gray and Hervella's
ones and ${\cal W}_5 \cong T^*M$. Besides the additional summand
${\cal W}_5$,  another interesting difference may be pointed out:
all the information about the torsion of the ${\sl
SU}(n)$-structure, $n \geq 4 $, is contained in the exterior
derivatives $d \omega$ and $d \psi_+$, or $d \omega$ and $d
\psi_-$. This happens similarly for another $G$-structures, $d
\varphi$ is sufficient to classify a ${\sl Spin}(7)$-structure, $d
\Omega$ is sufficient to know the intrinsic ${\sl Sp}(n){\sl
Sp}(1)$-torsion, $n > 2$, etc. However, we recall that $d \omega$
is not enough to classify a ${\sl U}(n)$-structure, we also need
to search in the Nijenhuis tensor for the remaining information.
Moreover, the importance of ${\sl SU}(n)$-structures from the
point of view of  geometry and theoretical physics  makes valuable
a detailed description of the involved tensors $\nabla \omega$ and
$\nabla \Psi$. Here we describe $\nabla \Psi$ which complements
the  study of $\nabla \omega$ done by Gray and Hervella.

The paper is organised as follows. In Section 2, we start
discussing  basic results. Then we pay attention to the study of
special almost hermitian $2n$-manifolds of high dimensions,
$2n\geq 8$. However, some results involving the cases $n =2,3$ are
also given. For instance, for $n\geq2$, we prove the invariance
under conformal changes of metric of a certain one-form  related
with parts ${\cal W}_4$ and ${\cal W}_5$ of the intrinsic torsion.
This is a generalization of a Chiossi and Salamon's result for
${\sl SU}(3)$-structures \cite{ChiossiSalamon}.

In Section 3, we study special almost Hermitian manifold of low
dimensions. Such manifolds of six dimensions have been studied in
\cite{ChiossiSalamon}. Here we show some additional detailed
information.  When $n=1,2,3$, the number of special peculiarities
that occur is big enough to justify a separated exposition. In
particular, we prove that, for these manifolds, $d \omega$, $d
\psi_+$ {\it and} $d \psi_-$ are sufficient to know the intrinsic
torsion.

Finally, as examples  of $\SUn(2n)$-structures, we consider almost
hyperhermitian manifolds in Section 4. We show that the exterior
derivatives $d \omega_I$, $d \omega_J$ and $d \omega_K$ of the
K{\"a}hler forms are enough to compute the covariant derivatives
$\nabla \omega_I$, $\nabla \omega_J$ and $\nabla \omega_K$. This
implies Hitchin's result \cite{Hitchin:Riemann-surface} that if
$\omega_I$, $\omega_J$ and $\omega_K$ are closed, then they are
covariant constant, i.e., the manifold is hyperk{\"a}hler.
Furthermore, we prove that locally conformal hyperk{\"a}hler manifolds
are equipped with three $\SUn(2n)$-structures of type ${\cal W}_4
\oplus {\cal W}_5$, respectively  associated  with the almost
complex structures $I$, $J$ and $K$. As a consequence of this
result, we obtain an alternative proof of the Ricci flatness of
the metric of hyperk{\"a}hler manifolds.

{\it Acknowledgements.} This work is supported by a grant from MEC 
(Spain), project MTM2004-2644.

\seccion{Special almost Hermitian manifolds}{\indent}
\setcounter{proposition}{0} \setcounter{proposition}{0} An {\it
almost Hermitian} manifold is a $2n$-dimensional manifold $M$, $n
>0$, with a ${\sl U}(n)$-structure. This means that $M$ is equipped
with a Riemannian metric $\langle \cdot , \cdot \rangle$  and an
orthogonal almost complex structure $I$. Each fibre $T_m M$ of the
tangent bundle  can be consider as complex vector space by
defining $i x = Ix$.  We will write $T_m M_{\Bbb C}$ when we are
regarding $T_m M$ as such a space.

We define a Hermitian scalar product $\langle \cdot , \cdot
\rangle_{\Bbb C} = \langle \cdot , \cdot  \rangle + i
\omega(\cdot,\cdot)$, where $\omega$ is the K{\"a}hler form given by
$\omega (x, y) = \langle x, Iy \rangle$. The real tangent bundle
$TM$ is identified with the cotangent bundle $T^*M$ by the map $x
\to \langle \cdot , x \rangle=x$. Analogously, the conjugate
complex vector space $\overline{T_m M}_{\Bbb C}$ is identified
with the dual complex space $T^{\ast}_m M_{\Bbb C}$ by the map $x
\to \langle \cdot , x \rangle_{\Bbb C} = x_{\Bbb C}$. It follows
immediately that $x_{\Bbb C} = x + i Ix$.

If we consider the spaces $\Lambda^p T^*_m M_{\Bbb C}$ of
skew-symmetric complex forms, one can check  $x_{\Bbb C} \wedge
y_{\Bbb C} = (x + i Ix) \wedge (y + i Iy)$. There are natural
extensions of scalar products to $\Lambda^p T^*_m M$ and
$\Lambda^p T^*_m M_{\Bbb C}$, respectively defined by
\begin{eqnarray*}
\displaystyle \langle a,b\rangle & = &  \frac{1}{p!} \sum_{i_1,
\dots , i_p=1}^{2n}
a(e_{i_1}, \dots, e_{i_p}) b(e_{i_1},\dots, e_{i_p}), \\
\displaystyle  \langle a_{\Bbb C} ,b_{\Bbb C} \rangle_{\Bbb C} & =
& \frac{1}{p!} \sum_{i_1, \dots , i_p=1}^{n} a_{\Bbb C} (u_{i_1},
\dots, u_{i_p}) \overline{ b_{\Bbb C}(u_{i_1},\dots, u_{i_p})},
\end{eqnarray*}
where $e_1, \dots , e_{2n}$ is an orthonormal basis for real
vectors and $u_1, \dots, u_n$ is a unitary basis for complex
vectors.

 The following conventions will be
used in this paper.  If $b$ is a $(0,s)$-tensor, we write
\begin{equation} \label{ecuacionesb}
  \begin{array}{l}
    I_{(i)}b(X_1, \dots, X_i, \dots , X_s) = - b(X_1, \dots , IX_i, \dots ,
    X_s),\\[2mm]
    I b(X_1,\dots,X_s) = (-1)^sb(IX_1,\dots,IX_s),\\[2mm]
{\it i}_{I}b   =  (I_{(1)} + \dots + I_{(s)})b, \\[2mm]
\displaystyle L(b) = \sum_{1 \leq i < j \leq s } I_{(i)} I_{(j)}
b, \qquad s \geq 2.
  \end{array}
\end{equation}

A {\it special almost Hermitian} manifold is a $2n$-dimensional
manifold $M$ with an ${\sl SU}(n)$-structure. This means that $(M,
\langle \cdot, \cdot \rangle, I)$ is an almost Hermitian manifold
equipped with a complex volume form $\Psi = \psi_+ + i \psi_-$
such that $\langle  \Psi , \Psi \rangle_{\Bbb C} = 1$. Note that
$I_{(i)} \psi_+ = \psi_-$.

 If  $e_1, \dots
, e_n $ is a unitary basis for complex vectors such that $  \Psi
(e_1,\dots , e_n)=1 , $ i.e., $ \psi_+ (e_1, \dots , e_n) = 1$ and
$\psi_- (e_1 , \dots , e_n)=0$,  then $e_1, \dots , e_n, Ie_1,
\dots , Ie_n$ is an orthonormal basis for real vectors {\it
adapted} to the ${\sl SU}(n)$-structure. Furthermore, if $A$ is a
matrix relating two adapted basis  of an ${\sl SU}(n)$-structure,
then $A \in {\sl SU}(n) \subseteq {\sl SO}(2n)$. On the other
hand, it is straightforward to check
$$
 \omega^n =  (-1)^{n(n+1)/2} n! \, e_1 \wedge \dots \wedge e_n \wedge Ie_1 \wedge
 \dots \wedge Ie_n,
$$
 where $\omega^n = \omega \wedge \stackrel{(n)}{\dots} \wedge  \omega$.

 If we fix the form $Vol$ such that $(-1)^{n(n+1)/2} n!\, Vol = \omega^n$ as real volume form, it
follows next lemma.
\begin{lemma} \label{volumenes}
Let $M$ be a special almost Hermitian $2n$-manifold, then
\begin{enumerate}
\item[{\rm (i)}] $\psi_+ \wedge \omega = \psi_- \wedge \omega = 0
$;
 \item[{\rm (ii)}] for $n$ odd,  we have
$ \psi_+ \wedge \psi_- = - (-1)^{n(n+1)/2}  2^{n-1} \, Vol$ and $
\psi_+ \wedge \psi_+ = \psi_- \wedge \psi_- = 0$;
 \item[{\rm (iii)}] for $n$ even, we have  $
\psi_+ \wedge \psi_+ = \psi_- \wedge \psi_- =  (-1)^{n(n+1)/2}
2^{n-1} \, Vol$ and $ \psi_+ \wedge \psi_- = 0$;
 \item[{\rm (iv)}] for $n \geq
2$ and $1 \leq i < j \leq n$, $I_{(i)} I_{(j)} \psi_+ = - \psi_+$
and $I_{(i)} I_{(j)} \psi_- = - \psi_-$; and
 \item[{\rm (v)}] $ x
\wedge \psi_+ =  Ix \wedge \psi_-= - (Ix \lrcorner \psi_+) \wedge
\omega$ and  $x \lrcorner \psi_+ = Ix \lrcorner \psi_-$, for all
vector $x$, where $\lrcorner$ denotes the interior product.
\end{enumerate}
\end{lemma}
\Dt All parts follow by a straightforward way, taking the
identities
\begin{eqnarray}
&& \psi_+  =  {\sl Re} \left( e_{1{\Bbb C}} \wedge \dots \wedge
e_{n{\Bbb C}}\right),  \label{psi+base} \qquad \psi_-  =  {\sl Im}
\left( e_{1{\Bbb C}} \wedge \dots \wedge
e_{n{\Bbb C}} \right), \label{psi-base} \\
&& \omega  =  \sum_{i=1}^{n} I e_i \wedge e_i,
\end{eqnarray}
into account, where $e_1, \dots , e_n, Ie_1, \dots , Ie_n$ is an
adapted basis to the ${\sl SU}(n)$-structure. Note that parts (ii)
and (iii) can be given together by the equation  \vspace{1mm}

 $\;$ \hspace{3cm}
$n! \,  \Psi \wedge  \overline \Psi = i^n \,(-1)^{n(n-1)/2} \, 2^n
\, \omega^n$. \cqd

\vspace{3mm}

We will also need to consider the contraction of a $p$-form $b$ by
a skew-symmetric contravariant two-vector $x\wedge y$, i.e., $
(x\wedge y) \lrcorner b(x_1, \dots , x_{p-2} ) = b (x,y,x_1, \dots
, x_{p-2})$.  When $n\geq 2$, it is obvious that $(Ix \wedge y)
\lrcorner \psi_+ = - (x\wedge y) \lrcorner \psi_-$. Furthermore,
let us note that there are two Hodge star operators defined on
$M$. Such operators, denoted by $\ast$ and $\ast_{\Bbb C}$, are
respectively associated with the volume forms $Vol$ and $\Psi$.
\vspace{2mm}

Relative to the real Hodge star operator, we have the following
results.
\begin{lemma} \label{estrellas}
 For any one-form $\mu$ we have
$$
\begin{array}{l}
\ast \left( \ast ( \mu \wedge \psi_+ ) \wedge \psi_+ \right) =
\ast \left( \ast ( \mu \wedge \psi_- ) \wedge \psi_- \right)  = -2^{n-2} \mu, \\[2mm]
 \ast \left( \ast ( \mu \wedge \psi_- ) \wedge
\psi_+ \right) = - \ast \left( \ast ( \mu \wedge \psi_+ ) \wedge
\psi_- \right) =  2^{n-2} I\mu.
\end{array}
$$
\end{lemma}
\Dt The identities follow by direct computation, taking equations
(\ref{psi+base})  into account. \cqd \vspace{3mm}

We are dealing with $G$-structures where $G$ is a subgroup of the
linear group ${\sl GL}(m , {\Bbb R})$. If $M$ possesses a
$G$-structure, then there always exists a $G$-connection defined
on $M$. Moreover, if $(M^m ,\langle \cdot , \cdot \rangle)$ is an
orientable $m$-dimensional Riemannian manifold and  $G$ is a
closed and connected subgroup of ${\sl SO}(m)$,  then there exists
a unique metric $G$-connection $\widetilde{\nabla}$ such that
$\xi_x = \widetilde{\nabla}_x - \nabla_x$ takes its values in
${\frak g}^{\perp}$,  where ${\frak g}^{\perp}$ denotes the
orthogonal complement in ${\frak s}{\frak o}(m)$ of the Lie
algebra ${\frak g}$ of $G$ and  $\nabla$ denotes the Levi-Civita
connection \cite{Salamon:holonomy,Cleyton:thesis}. The tensor
$\xi$ is the {\it intrinsic torsion}  of the $G$-structure and
$\widetilde{\nabla}$ is called the {\it minimal $G$-connection}.

 For ${\sl
U}(n)$-structures, the minimal ${\sl U}(n)$-connection is given by
$ \widetilde{\nabla} = \nabla + \xi$, with
\begin{equation} \label{torsion:xi}
\xi_X Y = - \frac{1}{2} I \left( \nabla_X I \right) Y.
\end{equation}
see \cite{Falcitelli-Farinola-Salamon}. Since ${\sl U}(n)$
stabilises the K{\"a}hler form $\omega$, it follows that
$\widetilde{\nabla} \omega = 0$. Moreover, the equation $\xi_X (
IY) + I ( \xi_X Y) =0$ implies  $\nabla \omega = - \xi \omega \in
T^* M \otimes {\frak u}(n)^{\perp}$. Thus, one can identify the
${\sl U}(n)$-components of $\xi$ with the ${\sl U}(n)$-components
of $\nabla \omega$. \vspace{2mm}

For ${\sl SU}(n)$-structures, we have the decomposition ${\frak
s}{\frak o}(2n)= {\frak s}{\frak u}(n) + {\Bbb R} + {\frak
u}(n)^{\perp}$, i.e., ${\frak s}{\frak u}(n)^{\perp} = {\Bbb R} +
{\frak u}(n)^{\perp}$. Therefore, the intrinsic ${\sl
SU}(n)$-torsion  $\eta + \xi$ is such that $\eta \in T^* M \otimes
{\Bbb R} \cong T^* M$ and $\xi$ is still determined by equation
\ref{torsion:xi}. The tensors $\omega$, $\psi_+$ and $\psi_-$ are
stabilised by the ${\sl SU}(n)$-action, and $\overline{\nabla}
\omega = 0$, $\overline{\nabla} \psi_+ = 0$ and $\overline{\nabla}
\psi_- = 0$, where $\overline{\nabla} = \nabla + \eta + \xi$ is
the minimal ${\sl SU}(n)$-connection. Since $\overline{\nabla}$ is
metric and $\eta \in T^* M \otimes {\Bbb R}$, we have $ \langle Y
, \eta_X Z \rangle = (I \eta)(X) \omega(Y,Z)$, where $\eta$ on the
right side is  a one-form. Hence
 \begin{equation} \label{torsion:eta} \eta_X Y = I
\eta (X) IY.
\end{equation}

We can check  $\eta \omega = 0$, then  from $\overline{\nabla}
\omega = 0$ we  obtain:
\begin{enumerate}
\item[{\rm (i)}] for $n=1$, $ \nabla \omega = - \xi \omega \in T^*
M \otimes {\frak u}(1)^{\perp} = \{ 0 \}$;
  \item[{\rm (ii)}] for $n=2$, $ \nabla \omega = - \xi \omega \in T^* M \otimes {\frak
u}(2)^{\perp} = {\cal W}_2 + {\cal W}_4$;
 \item[{\rm (iii)}] for $n \geq 3$,
$ \nabla \omega = - \xi \omega \in T^* M \otimes {\frak
u}(n)^{\perp} = {\cal W}_1 + {\cal W}_2 + {\cal W}_3 + {\cal W}_4,
$
\end{enumerate}
where the summands ${\cal W}_i$ are the irreducible ${\sl
U}(n)$-modules given by Gray and Hervella in \cite{GrayHervella}
and $+$ denotes direct sum. In general, these spaces ${\cal W}_i$
are also irreducible as ${\sl SU}(n)$-modules. The only exceptions
are ${\cal W}_1$ and ${\cal W}_2$ when $n=3$. In fact, for that
case, we have the following decompositions into irreducible ${\sl
SU}(3)$-components,
$$ {\cal W}_i = {\cal W}_i^+ + {\cal W}_i^-,
\quad i=1,2,
$$
where the space ${\cal W}_i^+$ (${\cal W}_i^-$) consists in those
tensors $a \in {\cal W}_i \subseteq T^*M \otimes \Lambda^2 T^* M$
such that the bilinear form $r(a)$, defined by $2 r(a) = \langle x
\lrcorner \psi_+ , y \lrcorner a \rangle$, is symmetric
(skew-symmetric). \vspace{3mm}

On the other hand, since $\overline{\nabla} \psi_+ = 0$ and
$\overline{\nabla} \psi_- = 0$, we have  $ \nabla \psi_+ = - \eta
\psi_+ - \xi \psi_+$ and $ \nabla \psi_- = - \eta \psi_- - \xi
\psi_-$. Therefore, from equations (\ref{torsion:xi}) and
(\ref{torsion:eta}) we obtain the following expressions
\begin{equation} \label{torsiones}
\begin{array}{lll}
\qquad - \eta_X \psi_+   =  - n I\eta(X) \psi_-, & \quad &  -\xi_X
\psi_+  = \displaystyle \frac{1}{2}  (e_i \lrcorner \nabla_X
\omega) \wedge
(e_i \lrcorner \psi_-), \\[2mm]
\qquad - \eta_X \psi_-   =    n I\eta(X) \psi_+,  & &
 -\xi_X \psi_-  =
 - \displaystyle \frac{1}{2}  (e_i \lrcorner \nabla_X \omega)
\wedge (e_i \lrcorner \psi_+),
\end{array}
\end{equation}
where the summation convention is used. \vspace{2mm}

 It is obvious that $- \eta \psi_+ \in {\cal W}^-_{5}
 = T^* M \otimes \psi_-$  and
 $- \eta \psi_- \in {\cal W}^+_{5} =  T^* M \otimes \psi_+$.
The tensors $-\xi \psi_+$ and $-\xi \psi_-$ are described in the
following proposition, where we need to consider the two ${\sl
SU}(n)$-maps
$$
\Xi_+
, \Xi_- \, : \, T^* M \otimes {\frak u}(n)^{\perp} \to T^* M
\otimes \Lambda^n T^* M
$$
 respectively defined  by $\nabla_{\cdot} \omega \to 1/2 \, (e_i \lrcorner
\nabla_{\cdot} \omega) \wedge (e_i \lrcorner \psi_-)$
 and
 $ \nabla_{\cdot} \omega \to - 1/2 \,
 (e_i \lrcorner \nabla_{\cdot} \omega) \wedge (e_i \lrcorner \psi_+)$.
Likewise, we also  consider the $SU(n)$-spaces $\lcf \lambda^{p,0}
\rcf  = \{ Re\left( b_{\Bbb C} \right) \, | \, b_{\Bbb C} \in
\Lambda^{p} T^* M_{\Bbb C} \}$ of real $p$-forms. Thus, $\lcf
\lambda^{0,0} \rcf = {\Bbb R}$, $\;\lcf \lambda^{1,0} \rcf = T^*
M$ and, for $p \geq 2$, $\lcf \lambda^{p,0} \rcf  = \{ b \in
\Lambda^{p} T^* M \, | \, I_{(i)}I_{(j)} b = - b, \, 1 \leq i < j
\leq p \}$. We write $\lcf \lambda^{p,0} \rcf$ in agreeing with
notations used in
\cite{Salamon:holonomy,Falcitelli-Farinola-Salamon}.
\begin{proposition} \label{Ximasmenos}
For $n \geq 3$, the $SU(n)$-maps $\Xi_+$ and $\Xi_-$ are injective
and
$$
\Xi_{+} \left( T^* M \otimes {\frak u}(n)^{\perp} \right) =
\Xi_{-} \left( T^* M \otimes {\frak u}(n)^{\perp} \right) = T^* M
\otimes \lcf \lambda^{n-2,0} \rcf  \wedge \omega.
$$
For $n=2$,  the maps $\Xi_+$ and $\Xi_-$ are not injective, being
$$
\begin{array}{l}
{\sl ker} \, \Xi_+ = T^* M \otimes \psi_-, \qquad {\sl ker} \,
\Xi_- = T^* M \otimes \psi_+, \\[3mm]
\Xi_{+} \left( T^* M \otimes {\frak u}(2)^{\perp} \right) =
\Xi_{-} \left( T^* M \otimes {\frak u}(2)^{\perp} \right) = T^* M
\otimes  \omega.
\end{array}
$$
\end{proposition}
\Dt We consider $n \geq 2$. As the real metric $\langle \cdot ,
\cdot \rangle$ is Hermitian with respect to $I$, we have $I (
\nabla_X \omega ) = - \nabla_X \omega$ \cite{GrayHervella}, for
all vector $X$. But this is equivalent to
$$
\nabla_X \omega = \sum^{n}_{1 \leq i < j \leq n} \left( a_{ij}
Re(e_{i{\Bbb C}} \wedge e_{j{\Bbb C}})   +  b_{ij} Im(e_{i{\Bbb
C}} \wedge e_{j{\Bbb C}}) \right)  \in \lcf \lambda^{2,0} \rcf ,
$$
where $e_1, \dots, e_n, Ie_1, \dots , Ie_n$ is an adapted basis.
Taking (\ref{torsiones}) into account, it is straightforward to
check
$$
\begin{array}{l}
\displaystyle \Xi_{+} \left( \nabla_X \omega \right) = -
\sum^{n}_{1 \leq i < j \leq n} a_{ij} Re\left( \ast_{\Bbb C}
\left( e_{i{\Bbb C}} \wedge
e_{j{\Bbb C}}\right) \right) \wedge \omega + \\
\displaystyle \hspace{1.85cm} + \sum^{n}_{1 \leq i < j \leq n}
b_{ij} Im\left( \ast_{\Bbb C} \left( e_{i{\Bbb C}} \wedge
e_{j{\Bbb C}}\right) \right) \wedge \omega \in \lcf \lambda^{2,0}
\rcf \wedge \omega, \\
\displaystyle \Xi_{-} \left( \nabla_X \omega \right) = -
\sum^{n}_{1 \leq i < j \leq n} a_{ij} Im\left( \ast_{\Bbb C}
\left( e_{i{\Bbb C}} \wedge
e_{j{\Bbb C}}\right) \right) \wedge \omega - \\
\displaystyle \hspace{1.85cm} - \sum^{n}_{1 \leq i < j \leq n}
b_{ij} Re\left( \ast_{\Bbb C} \left( e_{i{\Bbb C}} \wedge
e_{j{\Bbb C}}\right) \right) \wedge \omega \in \lcf \lambda^{2,0}
\rcf \wedge \omega.
\end{array}
$$
From  these equations  Proposition follows. \cqd \vspace{4mm}

For sake of simplicity, for $n \geq 2$, we denote ${\cal W}^{\Xi}=
T^* M \otimes \lcf \lambda^{n-2,0} \rcf \wedge \omega$. Moreover,
we  will consider the map $ {\cal L} \, : T^* M \otimes \Lambda^n
T^* M \to T^* M \otimes \Lambda^n T^* M $ defined by
\begin{equation} \label{laele}
{\cal L}(b) =  I_{(1)} \left(I_{(2)} + \dots + I_{(n+1)} \right)
b.
\end{equation}
 Proposition \ref{Ximasmenos} and above considerations give rise to the  following
theorem  where  we describe the properties satisfied by   the 
${\sl SU}(n)$-components of $\nabla \psi_+$ and $\nabla \psi_-$.
\begin{theorem} \label{mayor4}
Let $M$ be a special almost Hermitian $2n$-manifold, $n \geq 4$,
with K{\"a}hler form $\omega$ and complex volume form $\Psi = \psi_+ +
i \psi_-$. Then
$$
\begin{array}{l}
\nabla \psi_+ \in {\cal W}^{\Xi}_1 +  {\cal W}^{\Xi}_2 + {\cal
W}^{\Xi}_3 + {\cal W}^{\Xi}_4 + {\cal W}^-_5, \\[2mm]
 \nabla \psi_-
\in {\cal W}^{\Xi}_1 + {\cal W}^{\Xi}_2 + {\cal W}^{\Xi}_3 + {\cal
W}^{\Xi}_4 + {\cal W}^+_5,
\end{array}
$$
where ${\cal W}^{\Xi}_i = \Xi_{+} ({\cal W}_i) = \Xi_{-} ({\cal
W}_i)$, ${\cal W}_{5}^+ = T^*M \otimes \psi_{+}$ and ${\cal
W}_{5}^- = T^*M \otimes \psi_{-}$. The modules ${\cal W}^{\Xi}_i$
are explicitly described by
$$
\begin{array}{l}
{\cal W}^{\Xi}_1 = \{ e_i \otimes Ie_i \wedge b \wedge \omega +
e_i \otimes e_i \wedge I_{(1)} b \wedge \omega  \, | \, b
\in \lcf \lambda^{n-3,0} \rcf \}, \\[2mm]
{\cal W}^{\Xi}_2  = \{ b \in {\cal W}^{\Xi} \, | \, {\cal L}(b) =
(n-2) b \, \mbox{ and } \, \, \widetilde{\it a}(b) \wedge \omega = 0 \}, \\[2mm]
{\cal W}^{\Xi}_1 + {\cal W}^{\Xi}_2 = \{ b \in {\cal W}^{\Xi} \, |
\, {\cal L}(b) = (n-2) b \}, \\[2mm]
 {\cal
W}^{\Xi}_3  = \{ b \in {\cal W}^{\Xi} \, | \, \widetilde{\it a}(b) = 0 \}, \\[2mm]
{\cal W}^{\Xi}_4  = \{ e_i \otimes \left( (x \wedge e_i )
\lrcorner \psi_+ \right) \wedge \omega \, | \, x \in TM \} =\{ e_i
\otimes \left( (x \wedge e_i ) \lrcorner \psi_- \right) \wedge
\omega  \, | \, x \in TM \}, \\[2mm]
{\cal W}^{\Xi}_3 + {\cal W}^{\Xi}_4 = \{ b \in {\cal W}^{\Xi} \, |
\, {\cal L}(b) = - (n-2) b \},
\end{array}
$$
where $\widetilde{\it a}$ denotes the alternation map.
\end{theorem}
\Dt Some  parts of  Theorem follow  by computing the image $\Xi_+
\left( \nabla \omega \right)_{i}$ of the ${\cal W}_i$-part of
$\nabla \omega$, taking the properties for  ${\cal W}_i$ given in
\cite{GrayHervella} into account, and others, with Schur's Lemma
\cite{BtD} in mind, by computing $\Xi_{+} (a)$, where $0 \neq a
\in {\cal W}_i$. \cqd \vspace{2mm}

If we consider the alternation maps $\widetilde{\it a}_{\pm} \, : \,
{\cal W}^{\Xi}+ {\cal W}_5^{\mp} \to \Lambda^{n+1} T^* M$, we get
the following consequences of Theorem \ref{mayor4}.
\begin{corollary} \label{diferencialespsi}
For $n \geq 4$, the exterior derivatives of $\psi_+$ and
$\psi_{-}$ are such that
$$
d \psi_+, d \psi_- \in T^* M \wedge \lcf \lambda^{n-2,0} \rcf
\wedge \omega = {\cal W}_1^{a} + {\cal W}_2^{a} + {\cal
W}_{4,5}^{a},
$$
where  $\widetilde{\it a}_{\pm}({\cal W}_1^{\Xi}) ={\cal W}_1^{a}$,
$\widetilde{\it a}_{\pm}({\cal W}_2^{\Xi}) ={\cal W}_2^{a}$ and
$\widetilde{\it a}_{\pm}({\cal W}_4^{\Xi}) = \widetilde{\it a}_{\pm}({\cal
W}_5^{\mp}) = {\cal W}_{4,5}^{a}$. Moreover, ${\sl
ker}(\widetilde{\it a}_{\pm}) = {\cal W}^{\Xi}_3 + {\cal A}_{\pm}$, where
$T^* M \cong {\cal A}_{\pm} \subseteq {\cal W}_4^{\Xi} + {\cal
W}_5^{\mp}$, and the modules ${\cal W}_i^{a}$ are described by
$$
\begin{array}{l} {\cal W}_1^{a} = \lcf \lambda^{n-3,0} \rcf \wedge \omega \wedge \omega, \\[2mm]
{\cal W}_2^{a}   = \{ b \in T^*M \wedge \lcf \lambda^{n-2,0} \rcf
\wedge \omega \, | \, b \wedge \omega = 0 \,
\mbox{and} \, \ast b \wedge \psi_+ = 0 \, \}\\[2mm]
 \qquad = \{ b \in T^*M \wedge \lcf \lambda^{n-2,0} \rcf \wedge \omega \, | \, b \wedge \omega = 0  \,
\mbox{and} \, \ast b \wedge \psi_- = 0 \, \}, \\[2mm]
 {\cal W}_{4,5}^{a} = T^*M
\wedge \psi_+ = T^* M \wedge \psi_- = \lcf \lambda^{n-1,0} \rcf
\wedge \omega .
\end{array}
$$
\end{corollary}
Note also that
$$
\begin{array}{l}
{\cal W}_1^{a} + {\cal W}_2^{a} = \{ b \in T^*M \wedge \lcf
\lambda^{n-2,0} \rcf \wedge \omega \; | \, \ast
b \wedge \psi_+ = 0 \} \\[2mm]
\; \qquad \qquad \,  \, = \{ b \in T^*M \wedge \lcf \lambda^{n-2,0} \rcf
 \wedge \omega \; | \, \ast b \wedge \psi_- = 0 \},\\[2mm]
{\cal W}_2^{a} + {\cal W}_{4,5}^{a} = \{ b \in T^*M \wedge \lcf
\lambda^{n-2,0} \rcf \wedge \omega \, | \, b \wedge \omega = 0 \}.
\end{array}
$$

In this point we already have all the ingredients to explicitly
describe the one-form $\eta$. This will complete the definition of
 the ${\sl SU}(n)$-connection $\overline{\nabla}$.
\begin{theorem} \label{torsionw5}
For an ${\sl SU}(n)$-structure, $n \geq 2$,
 the ${\cal W}_5$-part $\eta$ of the torsion can be
identified with $- \eta \psi_+ = - n I \eta \otimes \psi_-$ or $-
\eta \psi_- =  n I \eta \otimes \psi_+$,
 where $\eta$ is a one-form such that
$$
\ast \left( \ast  d \psi_+ \wedge  \psi_+  +  \ast d \psi_- \wedge
\psi_- \right) =  n 2^{n-1} \eta + 2^{n-2} I d^* \omega,
$$
or
$$
\ast \left( \ast  d \psi_+ \wedge  \psi_- -  \ast d \psi_- \wedge
\psi_+ \right) =  n 2^{n-1} I \eta - 2^{n-2} d^* \omega .
$$
Furthermore, if $n\geq 3$, then $  \ast  d \psi_+ \wedge \psi_+ =
 \ast  d \psi_- \wedge \psi_- $ and $ \ast d
\psi_+ \wedge  \psi_-  \linebreak = -  \ast  d \psi_- \wedge
\psi_+ $.
\end{theorem}
\Dt We prove the result for $n\geq 4$ and we will see the cases
$n=2,3$ in next section. The ${\cal W}_4$-part of $\nabla \omega$
 is given by $ 2 (n-1) \left( \nabla \omega
\right)_{4} = e_i \otimes e_i \wedge d^* \omega  + e_i \otimes I
e_i \wedge I d^* \omega$ \cite{GrayHervella}. Then, by computing $
\Xi_{+} \left( \nabla \omega \right)_{4}$, we get
\begin{equation} \label{torsionw4}
 \left( \nabla \psi_+
\right)_{4} =  - \frac{1}{2(n-1)} e_i \otimes \left( (d^* \omega
\wedge e_i) \lrcorner \psi_+ \right) \wedge \omega.
\end{equation}
Now, since $\left( \nabla \psi_+ \right)_{5} = - n I \eta \otimes
\psi_-$, we have
$$
\widetilde{a}_+ \left( \left( \nabla \psi_+ \right)_{4} +  \left(
\nabla \psi_+ \right)_{5} \right)   =  - \frac{1}{2} \left( d^*
\omega \lrcorner \psi_+ \right) \wedge \omega - n I \eta \wedge
\psi_{-}
 =  - \left( \frac{1}{2} I d^* \omega + n \eta  \right) \wedge
\psi_+ .
$$
Hence, the ${\cal W}_{4,5}^{a}$-part of $d \psi_{+}$ is given by
\begin{equation} \label{nucleo}
 \left( d  \psi_+ \right)_{4,5}  = - \left(
n \eta + \frac{1}{2} I d^* \omega \right)  \wedge \psi_{+},
\end{equation}
Finally, taking Lemma \ref{estrellas} into account, it follows
$$
\begin{array}{l}
\ast \left( \ast  d \psi_+ \wedge  \psi_+ \right)  = \ast \left(
\left( d \psi_+ \right)_{4,5} \wedge \psi_+
\right) = n 2^{n-2} \eta + 2^{n-3} I d^* \omega , \\[2mm]
 \ast \left( \ast  d \psi_+ \wedge  \psi_- \right)  = \ast \left(
\left( d \psi_+ \right)_{4,5} \wedge \psi_- \right) =  n 2^{n-2} I
\eta - 2^{n-3} d^* \omega .
\end{array}
$$
The identities for $d \psi_{-}$ can be proved in a similar way.
\cqd \vspace{3mm}

\begin{observation}
\begin{enumerate}
\item[{\rm (i)}] {\rm It is known that $ I d^* \omega =  \ast (
\ast d \omega \wedge \omega) = - \langle \cdot \lrcorner d \omega
, \omega \rangle$. Therefore, Theorem \ref{torsionw5} says that,
for $n \geq 3$,  $\eta$ can be computed from $d \omega$ and $d
\psi_+$ ( or $d \psi_-$). For $n=2$, we will need $d\omega$, $d
\psi_+$ {\it and}  $d \psi_-$ to determine the one-form $\eta$.
\item[{\rm (ii)}] From equation (\ref{nucleo}), it follows that
${\cal A}_+ \subseteq {\sl ker} (\widetilde{a}_+)$ is given by
$$
{\cal A}_+ = \left\{  - \frac{1}{2(n-1)} e_i \otimes \left( (x
\wedge e_i) \lrcorner \psi_+ \right) \wedge \omega  - \frac{1}{2}
x \otimes \psi_- \, | \, x \in TM \right\}.
$$
Analogously, for ${\cal A}_- \subseteq {\sl ker}
(\widetilde{a}_-)$, we have
$$
{\cal A}_- =  \left\{ - \frac{1}{2(n-1)} e_i \otimes \left( (x
\wedge e_i) \lrcorner \psi_- \right) \wedge \omega  + \frac{1}{2}
x \otimes \psi_+ \, | \, x \in TM \right\}.
$$ }
\end{enumerate}
\end{observation}
\vspace{2mm}

 Since $ d \omega \in {\cal W}_1 + {\cal W}_3 + {\cal W}_4$
 and $d \psi_+, d \psi_- \in {\cal W}_1^{a} + {\cal W}_2^{a} + {\cal W}_{4,5}^{a}$,
 all the information about the
intrinsic torsion of an ${\sl SU}(n)$-structure, $n \geq 4$, is
contained in $d \omega$ and $d \psi_+$ (or $d \psi_-$). We recall
that, for a ${\sl U}(n)$-structure, $n \geq 2$, we need the
Nijenhuis tensor and $d \omega$ to have the complete information
about the intrinsic torsion. Equation (\ref{torsionw4}) and
Theorem \ref{torsionw5} give us the components ${\cal W}_4$ and
${\cal W}_5$ of $\nabla \psi_+$  in terms of  $d \omega$ and $d
\psi_+$. For sake of completeness, we will compute the remaining
parts of $\nabla \psi_+$ in terms of $d \omega$ and $d \psi_+$. To
achieve this, let us study the behavior of  the coderivatives $d^*
\psi_+$, $d^* \psi_-$ and the forms $d^*_{\omega} \psi_+$ and
$d^*_{\omega} \psi_-$  respectively defined by the contraction  of
$\nabla \psi_+$ and $\nabla \psi_-$ by $\omega$, i.e.,
$$
d^*_{\omega} \psi_+ (Y_1, \dots , Y_{n-1} ) = \nabla_{e_i} \psi_+
(Ie_i, Y_1, \dots , Y_{n-1} )
$$
 and an analog expression gives $d^*_{\omega} \psi_-$.

  Note that $d^* \psi_+ = - \ast d \ast \psi_+$ and $d^*
\psi_- = - \ast d \ast \psi_-$. By Lemma \ref{volumenes},  when
$n$ is odd (even), $\ast \psi_+ =- (-1)^{n(n+1)/2} \psi_-$ and
$\ast \psi_- = (-1)^{n(n+1)/2} \psi_+$ ($ \ast \psi_+ =
(-1)^{n(n+1)/2} \psi_+$ and $\ast \psi_- = (-1)^{n(n+1)/2} \psi_-
$). Therefore, by Corollary \ref{diferencialespsi}, it is
immediate that
$$
d^* \psi_+, d^* \psi_- \in  \ast \left( T^* M \wedge \lcf
\lambda^{n-2,0} \rcf \wedge \omega \right) = {\cal W}_1^{c} +
{\cal W}_2^{c} + {\cal W}_{4,5}^{c},
$$
where the modules ${\cal W}^c_{i}$ are described in the following
lemma.
\begin{lemma} \label{modulcoderivada}
 For $n \geq 4$, ${\cal W}^{c}=\ast \left(
T^* M \wedge \lcf \lambda^{n-2,0} \rcf \wedge \omega \right)$ and
$L$ the map defined by (\ref{ecuacionesb}),  the modules ${\cal
W}_1^{c}$, ${\cal W}_2^{c}$ and ${\cal W}_{4,5}^{c}$ are defined
by:
$$
\begin{array}{l}
{\cal W}_1^{c} = \lcf \lambda^{n-3,0} \rcf \wedge \omega,
\\[2mm]
{\cal W}_2^{c} = \{ a \in {\cal W}^{c} \, | \, a \wedge \omega
\wedge \omega = 0 \, \, \mbox{ and } \, \,  a \wedge \psi_+
=0\}, \\[2mm]
 {\cal W}_1^{c} + {\cal W}_2^{c}   =
\{ a \in {\cal W}^{c} \, | \, -2 L(a) = (n-2) (n-5) a \} = \{ a
\in {\cal W}^{c}
 \, | \, a \wedge \psi_+ =0 \}, \\[2mm]
{\cal W}_{4,5}^{c}  = \lcf \lambda^{n-1,0} \rcf = \{ x \lrcorner
\psi_{+} \, | \, x \in TM
\}, \\[2mm]
{\cal W}_2^{c} + {\cal W}_{4,5}^{c} = \{ a \in {\cal W}^{c} \, |
\, a \wedge \omega \wedge \omega = 0 \}.
\end{array}
$$
\end{lemma}
\Dt It follows  by applying  $\ast$ to the ${\cal W}_i^{a}$
modules  of Corollary \ref{diferencialespsi}.

For the description of ${\cal W}_1^{c} + {\cal W}_2^{c}$ involving
the map $L$. Taking Proposition \ref{Ximasmenos} into account, we
consider  $ \nabla \psi_+ = x \otimes b \wedge \omega \in T^* M
\otimes \lcf \lambda^{n-2,0} \rcf \wedge \omega. $ Now, making use
of Theorem \ref{mayor4}, we obtain the ${\cal W}_1^{\Xi} + {\cal
W}_2^{\Xi}$-part of $\nabla \psi_+$,
\begin{eqnarray} \label{nabla12psi+}
 (\nabla \psi_+)_{1,2} & = & \frac{(n-2) \nabla \psi_+ + {\cal L}
(\nabla \psi_+)}{2(n-2)} \\
& = & \frac{1}{2} \left( x \otimes b \wedge \omega + I x \otimes
I_{(1)} b \wedge \omega \right). \nonumber
\end{eqnarray}
Then, we  compute $(d^* \psi_+)_{1,2} = d^* (\nabla \psi_+)_{1,2}
$ and check  $  - 2 L (d^* \psi_+)_{1,2} = \linebreak (n-2)(n-5)
(d^* \psi_+)_{1,2}$. \cqd \vspace{3mm}

 In the following result,  $(d^* \psi_+)_i$ and
$(d^*_{\omega} \psi_+)_i$  ( $(d^* \psi_-)_i$ and $(d^*_{\omega}
\psi_-)_i$ )  respectively denote the images by the maps $d^*$ and
$d^*_{\omega}$  of the ${\cal W}_i$-component of $\nabla \psi_+$
($\nabla \psi_-$). This notation for the ${\cal W}_i$-part of a
tensor will be used in the following.
\begin{lemma} \label{lemaidentidadescoder}
For $n \geq 3$ and the map $i_I$ given by (\ref{ecuacionesb}), the
forms $d^* \psi_+$, $d^* \psi_-$, $d^*_{\omega} \psi_+$ and
$d^*_{\omega} \psi_+$ satisfy:
$$
\begin{array}{ll}
\left( d^{\ast} \psi_+ \right)_{1,2} = - \left( d^{\ast}_{\omega}
\psi_- \right)_{1,2}, \qquad &\left( d^{\ast}_{\omega}  \psi_+
\right)_{1,2} = \left( d^{\ast} \psi_- \right)_{1,2}, \\[2mm]
\left( d^{\ast} \psi_+ \right)_{4} = \left( d^{\ast}_{\omega}
\psi_- \right)_{4}, \qquad &\left( d^{\ast}_{\omega}  \psi_+
\right)_{4} = - \left( d^{\ast} \psi_- \right)_{4}, \\[2mm]
\left( d^{\ast} \psi_+ \right)_{5} = - \left( d^{\ast}_{\omega}
\psi_- \right)_{5} = n  \eta \lrcorner \psi_+, \qquad & \left(
d^{\ast}_{\omega} \psi_+ \right)_{5} = - \left( d^{\ast} \psi_-
\right)_{5} =  n \eta \lrcorner \psi_-, \\[2mm]
i_I \left( d^{\ast} \psi_{\pm} \right)_{1,2} = (n-3) \left(
d^{\ast}_{\omega} \psi_{\pm} \right)_{1,2}, \qquad & i_I \left(
d^{\ast}_{\omega}  \psi_{\pm}
\right)_{1,2} = - (n-3) \left( d^{\ast} \psi_{\pm} \right)_{1,2}, \\[2mm]
i_I \left( d^{\ast} \psi_{\pm} \right)_{4} = - (n-1) \left(
d^{\ast}_{\omega} \psi_{\pm} \right)_{4}, \qquad & i_I \left(
d^{\ast}_{\omega}  \psi_{\pm} \right)_{4} = (n-1) \left( d^{\ast}
\psi_{\pm} \right)_{4}.
\end{array}
$$
\end{lemma}
\Dt Here we only consider  $n\geq 4$, the proof for $n=3$ will be
shown in next section. The identities of fourth and fifth lines
follow by similar arguments to those contained in the proof of
Lemma \ref{modulcoderivada}. The identities of the third line
follow by a straightforward way.

Making use of the maps $\Xi_+$, $\Xi_-$ and equations
(\ref{torsiones}), we note
\begin{equation} \label{coderpsi-via-psi+}
\begin{array}{l}
\left( \nabla_{\cdot} \psi_- \right)_{1,2} = \left( \nabla_{I \cdot }
\psi_+ \right)_{1,2}, \qquad \left( \nabla_{\cdot} \psi_- \right)_{3,4}
= - \left( \nabla_{I \cdot } \psi_+ \right)_{3,4}.
\end{array}
\end{equation}
 Hence, applying the maps $d^*$ and $d^*_{\omega}$ to
both sides of these equalities, the  identities of first and
second lines in Lemma follow. \cqd \vspace{2mm}

We know how to compute $(\nabla \psi_+)_4$ and $(\nabla \psi_+)_5$
(equation (\ref{torsionw4}) and Theorem \ref{torsionw5}) . Now, we
will  show expressions for the remaining ${\sl SU}(n)$-parts of
$\nabla \psi_+$ in terms of $d \omega$ and $d \psi_+$.
\begin{proposition} \label{calculoexolic}
Let $M$ be a special almost Hermitian $2n$-manifold, $n \geq 4$.
Then
\begin{enumerate}
\item[{\rm (i)}]$\left( \nabla_{ \cdot } \psi_+ \right)_{1} = e_i
\otimes Ie_i \wedge b \wedge \omega + e_i \otimes e_i \wedge
I_{(1)} b \wedge \omega$, $\left( d \psi_+ \right)_{1} = - 2 b
\wedge \omega \wedge \omega$, $\left( d^*  \psi_+ \right)_{1} = 2
I_{(1)} b \wedge \omega$ and $\; \left( d^*_{\omega}  \psi_+
\right)_{1} = - 2 b \wedge \omega$,
 where $b$ is
given by $b = Im\left( \ast_{\Bbb C} \left( \left( \nabla \omega
\right)_{1} + i I_{(1)} \left( \nabla \omega\right)_{1} \right)
\right)$, $\;\;12 (-1)^{n} b = I \ast ( d \psi_{+} \wedge \omega
)$, $\;\;12 (-1)^{n} I_{(1)} b =  I \ast ( d^* \psi_{+} \wedge
\omega \wedge \omega)$, or $12 (-1)^{n-1} b = I \ast (
d^*_{\omega} \psi_{+} \wedge \omega \wedge \omega )$;
 \item[{\rm (ii)}] $
\left( \nabla_{ \cdot } \psi_+ \right)_{2} =  e_i \otimes e_i
\lrcorner ( d^*_{\omega} \psi_+)_{1,2} \wedge \omega + e_i \otimes
I e_i \lrcorner (d^* \psi_+)_{1,2} \wedge \omega - 8 ( \nabla
\psi_+ )_1$, where \linebreak  $4 (n-2) (a)_{1,2} = (n-1)(n-2) a +
2 L(a)$, for $a=d^* \psi_+, d^*_{\omega} \psi_+$;
 \item[{\rm (iii)}] $ 2 \left( \nabla_{ \cdot } \psi_+
\right)_{3} =  \Xi_+ \left( ( 1 - I_{(2)} I_{(3)} )  (d
\omega)_{3} \right)$, where $(d \omega)_{3} = (d \omega)_{3,4} -(d
\omega)_{4}$ with $4 (d \omega)_{3,4} = 3 d \omega + L(d \omega)$
and $\; (n-1) (d \omega)_{4} = - I d^* \omega \wedge \omega$.
\end{enumerate}
\end{proposition}
\Dt By Theorem \ref{mayor4}, $ \left( \nabla_{ \cdot } \psi_+
\right)_{1} = e_i \otimes Ie_i \wedge b \wedge \omega + e_i
\otimes e_i \wedge I_{(1)} b \wedge \omega, $ where $b \in \lcf
\lambda^{n-3,0} \rcf$. Therefore, $ \left( d \psi_+ \right)_{1} =
- 2 b \wedge \omega \wedge \omega. $ On the other hand, it is not
hard to check
\begin{equation} \label{comeomeome}
\ast \left( c \wedge \omega \wedge \omega \wedge \omega \right) =
6 (-1)^{n-1} I c,
\end{equation}
for all $c \in \lcf \lambda^{n-3,0} \rcf$. Therefore, taking this
last  identity into account, we have
$$
\ast \left( d \psi_+ \wedge \omega \right)  = \ast \left( \left( d
\psi_+ \right)_{1} \wedge \omega \right) = 12 (-1)^{n-1} I b.
$$
Now, let us assume $ \nabla \omega = c \in \lcf \lambda^{3,0} \rcf
= {\cal W}_1$. Then, computing $\Xi_+ (c)$, we have
$$
\Xi_+ (c) = e_i \otimes Ie_i \wedge Im \left( \ast_{\Bbb C} \left(
c + i I_{(1)} c \right) \right) \wedge \omega  - e_i \otimes e_i
\wedge Re \left( \ast_{\Bbb C} \left( c + i I_{(1)} c \right)
\right) \wedge \omega.
$$
Hence the first identity for $b$ follows. The remaining identities
of (i)  involving $d^* \psi_+$ and $d^*_{\omega}\psi_+$ follow by
a straightforward way from $(\nabla \psi_+)_1$, taking equation
(\ref{comeomeome}) into account.

For part (ii).  If $ \nabla \psi_+ = x \otimes b \wedge \omega \in
T^*M \otimes \lcf \lambda^{n-2,0} \rcf$, by equation
(\ref{nabla12psi+}), we have $ 2 (\nabla \psi_+)_{1,2} = \left( x
\otimes b \wedge \omega + I x \otimes I_{(1)} b \wedge \omega
\right)$. Therefore, making use of part (i), it follows
$$
6 (\nabla \psi_+ )_1 = e_i \otimes I e_i \wedge (x \lrcorner
I_{(1)} b) \wedge \omega + e_i \otimes e_i \wedge I_{(1)} (x
\lrcorner I_{(1)} b) \wedge \omega.
$$
Moreover,
\begin{eqnarray} \label{compcoder12}
2 (d^* \psi_+)_{1,2} & = & I x \wedge b - x \wedge I_{(1)} b - 2
(x \lrcorner b) \wedge \omega, \\
 2 (d^*_{\omega} \psi_+)_{1,2} & = &  x
\wedge b + I x \wedge I_{(1)} b - 2 (x \lrcorner b) \wedge \omega.
\label{compcoderome12}
\end{eqnarray}
 From these equations, it is not hard to check
$$
   e_i \otimes e_i \lrcorner (d^*_{\omega} \psi_+)_{1,2} \wedge \omega + e_i \otimes I e_i
\lrcorner (d^* \psi_+)_{1,2} \wedge \omega = 2 (\nabla
\psi_+)_{1,2} + 6 ( \nabla \psi_+ )_1.
$$
Hence the first identity of (ii) follows. Furthermore, by Lemma
\ref{modulcoderivada}, we have  the equalities
$$
 - 2 L(d^* \psi_+)_{1,2}  =  (n-2)(n-5) (d^* \psi_+)_{1,2},
  \qquad - 2 L(d^* \psi_+)_{4,5}  =  (n-1)(n-2) (d^*
 \psi_+)_{4,5}.
$$
Therefore,
\begin{eqnarray*}
4 (n-2) (d^* \psi_+)_{1,2} & = & (n-1) (n-2) d^* \psi_+ + 2 L(d^*
\psi_+), \\
 4 (n-2) (d^* \psi_+)_{4,5} &  = & -(n-2) (n-5) d^* \psi_+ - 2 L(d^*
\psi_+).
\end{eqnarray*}
The required expression for $(d^*_{\omega} \psi_+)_{1,2}$ can be
deduced in a similar way.

Finally, part (iii) follows from identities for $\nabla \omega$
given in  \cite{Grayminimal} and \cite{GrayHervella}. \cqd

\begin{observation} {\rm
\begin{enumerate}
\item[(i)] From the identities given in Lemma
\ref{lemaidentidadescoder}, the forms $d^* \psi_+$ and
$d^*_{\omega} \psi_+$ can be computed in terms of $d \psi_+$ ($d
\psi_-$). Thus Proposition \ref{calculoexolic} corroborates our
claiming  that, for $n \geq 4$,  $d \omega$ and $d \psi_+$ ($d
\psi_-$) are enough to know the intrinsic ${\sl SU}(n)$-torsion.
 \item[(ii)]
 Taking equations (\ref{coderpsi-via-psi+}) into account,
it is not hard to deduce the respective ${\sl SU}(n)$-components,
$n\geq 4$, of $\nabla \psi_-$ from those of $\nabla \psi_+$.
\end{enumerate} }
\end{observation}

Relative with conformal changes of metric, we point out the
following facts which are  generalizations  of  results for ${\sl
SU}(3)$-structures proved by Chiossi and Salamon
\cite{ChiossiSalamon}.
\begin{proposition} \label{cambio4}
For conformal changes of metric given by $\langle \cdot , \cdot
\rangle_o = e^{2f} \langle \cdot , \cdot \rangle$, the ${\cal
W}_4$ and  ${\cal W}_5$ parts of the intrinsic ${\sl
SU}(n)$-torsion, $n \geq 2$, are modified in the way
$$
I d^* \omega_o = I d^* \omega - 2 (n-1) df, \qquad \eta_o = \eta -
\frac{1}{n} df,
$$
where $\omega_o$ and $\eta_o$ are respectively the K{\"a}hler form
and the ${\cal W}_5$ one-form of the metric $\langle \cdot , \cdot
\rangle_o$. Moreover, the one-form $ 2n(n-1) \eta -  I d^* \omega$
is not altered by such changes of metric.
\end{proposition}
\Dt  On one hand, the equation for $I d^* \omega_o$ was deduced in
\cite{GrayHervella}. On the other hand, from $\psi_{+o}= e^{nf}
\psi_+$ and  $\psi_{-o}= e^{nf} \psi_-$, we have $ d \psi_{+o} = n
e^{nf} d f \wedge \psi_+ + e^{nf} d \psi_+$ and $ d \psi_{-o} = n
e^{nf} d f \wedge \psi_- + e^{nf} d \psi_-$. Moreover, if $\ast_o$
is the Hodge star operator for $\langle \cdot , \cdot \rangle_o$
and $\alpha$ is a $p$-form, then $\ast_o \alpha = e^{2(n-p)f} \ast
\alpha$. Taking this last identity into account, we deduce
\begin{eqnarray*}
\ast_o \left( \ast_o d \psi_{+o} \wedge \psi_{+o} \right) + \ast_o
\left( \ast_o d \psi_{-o} \wedge \psi_{-o} \right) & = &  \ast
\left( \ast d \psi_+ \wedge \psi_+ \right)   + \ast \left( \ast d
\psi_- \wedge \psi_- \right)\\
&& - n 2^{n-1} d f.
\end{eqnarray*}
The required identity for $\eta_o$ follows from this last
identity and Theorem \ref{torsionw5}. Finally, it is obvious that
$ 2n(n-1) \eta_o - I d^* \omega_o = 2n(n-1) \eta - I d^* \omega. $
 \cqd

\begin{observation}{\rm
By Proposition \ref{cambio4}, for $n=3$,  the one-form $ 12 \eta -
I d^* \omega$ is not altered by conformal changes of metric. In
\cite{ChiossiSalamon}, Chiossi and
 Salamon consider six-dimensional manifolds with ${\sl
SU}(3)$-structure and prove that the tensor $3 \tau_{{\cal W}_4} +
2 \tau_{{\cal W}_5}$ is not modified under conformal changes of
metric, where $\tau_{{\cal W}_4}$ and $\tau_{{\cal W}_5}$ are
one-forms such that,  in the terminology here used, are given by $
2\tau_{{\cal W}_4} = -  I d^* \omega$ and  $2 \tau_{{\cal W}_5} =
 \eta  + I d^* \omega$.
Note that $ 3 \tau_{{\cal W}_4} + 2 \tau_{{\cal W}_5} =
\frac{1}{2}  \left( 12 \eta - I d^* \omega \right). $
 }
\end{observation}

\vspace{4mm}

\seccion{Low dimensions} \noindent \setcounter{proposition}{0} In
this section we consider special almost Hermitian manifolds of
dimension two, four and six. \vspace{4mm}

\noindent {\bf 3.1 Six dimensions}

\noindent  Here we  focus our attention on the very special case
of six-dimensional manifolds with  ${\sl SU}(3)$-structure (see
\cite{ChiossiSalamon}). In this case, we have
\begin{equation} \label{caso3}
\nabla \omega \in T^* M \otimes {\frak u}(3)^{\perp} = {\cal
W}^+_1 + {\cal W}^-_1 + {\cal W}_2^+ + {\cal W}_2^- + {\cal W}_3 +
{\cal W}_4.
 \end{equation}
 If we denote $[ T^* M_{\Bbb C} \otimes_{\Bbb C} \Lambda^2 T^*
M_{\Bbb C}] = \{ Re\left( b_{\Bbb C} \right) \,|\, b_{\Bbb C} \in
T^* M_{\Bbb C} \otimes_{\Bbb C} \Lambda^2 T^* M_{\Bbb C} \}$, some
summands in (\ref{caso3}) are described by
 \begin{eqnarray*}
 && {\cal W}^+_1 =
{\Bbb R} \psi_+, \qquad \;{\cal W}^-_1 = {\Bbb R} \psi_-, \\
&&{\cal W}_1^+ + {\cal W}^+_2 = \left\{ b \in [ T^* M_{\Bbb C}
\otimes_{\Bbb C} \Lambda^2 T^* M_{\Bbb C}] \,\, | \, \langle \cdot
\lrcorner \psi_+ , \cdot \lrcorner b \rangle \mbox{ is symmetric}
\right\},\\
&& {\cal W}_1^- + {\cal W}^-_2 = \left\{ b \in [ T^* M_{\Bbb C}
\otimes_{\Bbb C} \Lambda^2 T^* M_{\Bbb C}] \, \, | \, \langle
\cdot \lrcorner \psi_+ , \cdot \lrcorner b \rangle \mbox{ is
  skew-symmetric} \right\}.
\end{eqnarray*}

 By Proposition
\ref{Ximasmenos}, the ${\sl SU}(3)$-maps $\Xi_+$ and $\Xi_-$ are
injective and
$$
\Xi_+ \left( T^* M \otimes {\frak u}(3)^{\perp} \right) = \Xi_-
\left( T^* M \otimes {\frak u}(3)^{\perp} \right) = T^* M \otimes
T^* M \wedge \omega.
$$
In the following theorem we describe properties of the ${\sl
SU}(3)$-components of $\nabla \psi_{+}$ and $\nabla \psi_{-}$.
\begin{theorem}
 \label{igual3} Let $M$ be a special almost Hermitian
 $6$-manifold
 with K{\"a}hler form $\omega$ and complex volume form $\Psi =
\psi_+ + i \psi_-$. Then
$$
\begin{array}{l}
\nabla \psi_+ \in {\cal W}^{\Xi ; a}_1 + {\cal W}^{\Xi ; b}_1 +
{\cal W}^{\Xi ; a}_2 +  {\cal W}^{\Xi ; b}_2 + {\cal
W}^{\Xi}_3 + {\cal W}^{\Xi}_4 + {\cal W}^-_5, \\[2mm]
 \nabla \psi_-
\in {\cal W}^{\Xi ; a }_1 + {\cal W}^{\Xi ; b}_1  + {\cal W}^{\Xi
; a}_2 + {\cal W}^{\Xi ; b }_2 + {\cal W}^{\Xi}_3 + {\cal
W}^{\Xi}_4 + {\cal W}^+_5,
\end{array}
$$
where ${\cal W}^{\Xi;a}_i = \Xi_{+} ({\cal W}^+_i) = \Xi_{-}
({\cal W}^-_i)$, ${\cal W}^{\Xi;b}_i = \Xi_{+} ({\cal W}^-_i) =
\Xi_{-} ({\cal W}^+_i)$, $i=1,2$; ${\cal W}^{\Xi}_j = \Xi_{+}
({\cal W}_j) = \Xi_{-} ({\cal W}_j)$, $j=3,4$;
 ${\cal W}_{5}^+ = T^*M \otimes \psi_{+}$ and
${\cal W}_{5}^- = T^*M \otimes \psi_{-}$. If ${\cal W}^{\Xi} = T^*
M \otimes T^*M \wedge \omega$,  ${\cal L}$ is the map defined by
(\ref{laele}) and $\widetilde{a}$ denotes the alternation map, the
modules ${\cal W}^{\Xi;a}_i$, ${\cal W}^{\Xi;b}_i$ and ${\cal
W}^{\Xi}_j$ are described by
$$
\begin{array}{l}
{\cal W}^{\Xi;a}_{1} = {\Bbb R} e_i \otimes e_i \wedge \omega,
\qquad {\cal W}^{\Xi;b}_{1} = {\Bbb R} e_i \otimes Ie_i  \wedge
\omega,\\[2mm]
{\cal W}^{\Xi ;a}_{2} = \{  b \in  {\cal W}^{\Xi}  | \, \langle
b(e_i,e_i, \cdot , \cdot) , \omega \rangle = 0, \, b(e_i,Ie_i,
\cdot , \cdot) =0 \,  \mbox{ and } \, {\cal L}(b) =b
\, \}, \\[2mm]
{\cal W}^{\Xi ; b}_{2} = \{ b \in {\cal W}^{\Xi}  | \, \langle
b(e_i, Ie_i, \cdot , \cdot) , \omega  \rangle = 0, \, b(e_i, e_i,
\cdot , \cdot) = 0 \, \mbox{ and } \, {\cal L}(b)=b \},
\\[2mm]
{\cal W}^{\Xi;a}_1 + {\cal W}^{\Xi;b}_1 + {\cal W}^{\Xi;a}_2 +
{\cal W}^{\Xi;b}_2 = \{ b \in {\cal W}^{\Xi} \, | \, {\cal L}(b) =
b
\}, \\[2mm]
{\cal W}^{\Xi}_{3} = \{ b \in {\cal W}^{\Xi}  | \, {\cal L} (b) =
- b \, \mbox{ and } \; \widetilde{a}(b)=0
\},\\[2mm]
{\cal W}^{\Xi}_{4} = \{ e_i \otimes \left( (x \wedge e_i )
\lrcorner \psi_+ \right) \wedge \omega  |  x \in TM \} =\{ e_i
\otimes \left( (x \wedge e_i ) \lrcorner \psi_- \right) \wedge
\omega   | x \in TM \}, \\[2mm]
{\cal W}^{\Xi}_3 + {\cal W}^{\Xi}_4 = \{ b \in {\cal W}^{\Xi} \, |
\, {\cal L}(b) = -  b \}.
\end{array}
$$
\end{theorem}
\Dt We can proceed in a similar way as in the proof of Theorem
\ref{mayor4}. \cqd \vspace{2mm}

If we consider the alternation maps $\widetilde{a}_{\pm} \, : \, T^* M
\otimes T^* M \wedge \omega   + {\cal W}_5^{\mp} \to \Lambda^{4}
T^* M$, we get the following consequences of Theorem \ref{igual3}.
\begin{corollary}
For  ${\sl SU}(3)$-structures, the exterior derivatives of
$\psi_+$ and $\psi_-$ are such that
$$
d \psi_+, d \psi_- \in \Lambda^{4} T^* M = {\cal W}_1^{a} + {\cal
W}_2^{a} + {\cal W}_{4,5}^{a},
$$
where  $\widetilde{a}_{\pm} ({\cal W}_1^{\Xi;b}) = {\cal W}_1^{a}$,
$\widetilde{a}_{\pm}({\cal W}_2^{\Xi;b}) ={\cal W}_2^{a}$ and
$\widetilde{a}_{\pm}({\cal W}_4^{\Xi}) = \widetilde{a}_{\pm}({\cal
W}_5^{\mp}) = {\cal W}_{4,5}^{a}$. Moreover, ${\sl
Ker}(\widetilde{a}_{\pm}) = {\cal W}_1^{\Xi;a} + {\cal W}_2^{\Xi;a} +
{\cal W}^{\Xi}_3 + {\cal A}_{\pm}$, where $T^* M \cong {\cal A}_{\pm}
\subseteq {\cal W}_4^{\Xi} + {\cal W}_5^{\mp}$, and the modules
${\cal W}_i^{a}$ are described by
$$
\begin{array}{l} {\cal W}_1^{a} = {\Bbb R}  \omega \wedge \omega, \\[2mm]
{\cal W}_2^{a}   =  {\frak s}{\frak u}(3) \wedge \omega = \{ b \in
\Lambda^{4} T^* M \, | \, b \wedge
\omega = 0 \, \mbox{ and } \, \ast b \wedge \psi_+ = 0 \, \}\\[2mm]
 \qquad = \{ b \in \Lambda^{4} T^* M   \, | \, b \wedge \omega = 0  \,
\mbox{ and } \, \ast b \wedge \psi_- = 0 \, \}, \\[2mm]
 {\cal W}_{4,5}^{a} = T^*M
\wedge \psi_+ = T^* M \wedge \psi_- =  \lcf \lambda^{2,0} \rcf \wedge \omega \\[2mm]
\qquad  = \{ x \lrcorner \psi_+ \wedge \omega \, | \, x \in TM \}
= \{ x \lrcorner \psi_- \wedge \omega \, | \, x \in TM \}.
\end{array}
$$
\end{corollary}

Moreover,  we also have
$$
\begin{array}{l}
{\cal W}_1^{a} + {\cal W}_2^{a} = \{ b  \in \Lambda^{4}
T^* M  \; | \, \ast b \wedge \psi_+ = 0 \} = \{ b \in \Lambda^{4} T^* M  \; | \, \ast b \wedge \psi_- = 0 \},\\[2mm]
{\cal W}_2^{a} + {\cal W}_{4,5}^{a} = \{ b  \in \Lambda^{4} T^* M
\, | \, b \wedge \omega = 0 \}.
\end{array}
$$

In this point, one can proceed as  in the proof,  for high
dimensions,  of Theorem \ref{torsionw5} and obtain the results of
such Theorem for $n=3$. Along such a proof we  would get
\begin{eqnarray} \label{torsion3w4}
 \left( \nabla \psi_+
\right)_{4} & = &  \Xi_{+} \left( \nabla \omega \right)_{4} = -
\frac{1}{4} e_i \otimes \left( (d^* \omega \wedge e_i) \lrcorner
\psi_+ \right) \wedge \omega, \\
\label{nucleo3}
 \left( d  \psi_+ \right)_{4,5} & = &- \left(
3 \eta + \frac{1}{2} I d^* \omega \right)  \wedge \psi_{+}.
\end{eqnarray}
Likewise, in a similar way, we would also obtain
\begin{eqnarray}
 \left( \nabla \psi_-
\right)_{4} =  \Xi_{-} \left( \nabla \omega \right)_{4} & = & -
\frac{1}{4} e_i \otimes \left( (d^* \omega \wedge e_i) \lrcorner
\psi_- \right) \wedge \omega, \\ \label{nucleo3menos}
 \left( d  \psi_- \right)_{4,5} & = &  - \left(
3 \eta + \frac{1}{2} I d^* \omega \right)  \wedge \psi_{-}.
\end{eqnarray}

\begin{observation}
\begin{enumerate} {\rm
\item[{\rm (i)}] From equation (\ref{nucleo3}), it follows that
${\cal A}_+ \subseteq {\sl ker} (\widetilde{a}_+)$ is given by
$$
{\cal A}_+ = \left\{  - \frac{1}{4} e_i \otimes \left( (x \wedge
e_i) \lrcorner \psi_+ \right) \wedge \omega  - \frac{1}{2} x
\otimes \psi_- \, | \, x \in TM \right\}.
$$
Analogously, from equation (\ref{nucleo3}), for ${\cal A}_-
\subseteq {\sl ker} (\widetilde{a}_-)$, we have
$$
{\cal A}_- =  \left\{ - \frac{1}{4} e_i \otimes \left( (x \wedge
e_i) \lrcorner \psi_- \right) \wedge \omega  + \frac{1}{2} x
\otimes \psi_+ \, | \, x \in TM \right\}.
$$
 \item[{\rm (ii)}] Theorem \ref{torsionw5} says that  $\eta$ can
 be computed from $d \omega$ and $d \psi_+$ ($d \psi_-$).
Moreover, since $ d \omega \in {\cal W}^+_1 + {\cal W}^-_1 + {\cal
W}_3 + {\cal W}_4$  and
$$
\begin{array}{l}
d \psi_+  \in {\cal W}_1^{a} + {\cal W}_2^{a} + {\cal W}_{4,5}^{a}
=  \widetilde{a}_+ \com \Xi_+ \left( {\cal W}_1^- + {\cal W}_2^- +
{\cal W}_4 \right) + \widetilde{a}_+ \left( {\cal W}_5^-
\right), \\[2mm]
d \psi_-  \in {\cal W}_1^{a} + {\cal W}_2^{a} + {\cal W}_{4,5}^{a}
=  \widetilde{a}_- \com \Xi_- \left( {\cal W}_1^+ + {\cal W}_2^+ +
{\cal W}_4 \right) + \widetilde{a}_- \left( {\cal W}_5^+ \right),
\end{array}
$$
we need $d \omega$, $d \psi_+$ $and$ $d \psi_-$ to have the whole
information about the intrinsic  $\SUn(3)$-torsion. }
\end{enumerate}
\end{observation}

The  ${\cal W}_4$ and ${\cal W}_5$ parts of $\nabla \psi_+$ are
given by equation (\ref{torsion3w4}) and Theorem \ref{torsionw5}.
As in the previous section, for sake of completeness, we will see
how to compute the remaining parts of $\nabla \psi_+$ by using $d
\omega$, $d \psi_+$ $and$ $d \psi_-$. For such a purpose, we study
properties of the coderivatives $d^* \psi_+$, $d^* \psi_-$ and the
two-forms $d^*_{\omega} \psi_+$ and $d^*_{\omega} \psi_-$. Note
that,
 by Lemma \ref{volumenes},  we have $ d^* \psi_+ =  \ast d
 \psi_-$ and $ d^* \psi_- =   - \ast d
 \psi_+$. Therefore,
$$
d^* \psi_+, d^* \psi_-,d^*_{\omega} \psi_+, d^*_{\omega} \psi_-
\in \Lambda^2 T^* M = {\cal W}_1^{c} + {\cal W}_2^{c} + {\cal
W}_{4,5}^{c},
$$
where ${\cal W}_1^{c} = \ast \left( {\cal W}_1^{a}\right)$, ${\cal
W}_2^{c} = \ast \left( {\cal W}_2^{a}\right)$ and ${\cal
W}_{4,5}^{c} = \ast \left({\cal W}_{4,5}^{a}\right)$.
\begin{lemma} \label{modul3coderivada}
 For ${\sl SU}(3)$-structures,  the modules ${\cal
W}_1^{c}$, ${\cal W}_2^{c}$ and ${\cal W}_{4,5}^{c}$ are defined
by:
$$
\begin{array}{l}
{\cal W}_1^{c} = {\Bbb R} \omega, \quad {\cal W}_2^{c} = \{ b \in
\Lambda^2 T^* M \, | \, b \wedge \omega \wedge \omega = 0 \, \,
\mbox{and} \, \,  b \wedge \psi_+
=0\}, \\[2mm]
 {\cal W}_1^{c} + {\cal W}_2^{c}   =
\{ b \in \Lambda^2 T^* M \, | \, Ib = b \} = \{ b \in \Lambda^2 T^* M \, | \, b  \wedge \psi_+ =0 \}, \\[2mm]
{\cal W}_{4,5}^{c}  = \lcf \lambda^{2,0} \rcf = \{ x \lrcorner
\psi_{+} \, | \, x \in TM \},\\[2mm]
 {\cal W}_2^{c} + {\cal
W}_{4,5}^{c} = \{ b \in  \Lambda^2 T^* M \, | \, b \wedge \omega
\wedge \omega = 0 \}.
\end{array}
$$
\end{lemma}
\Dt It follows by similar arguments as in the proof of Lemma
\ref{modulcoderivada}. \cqd

Now one can prove the identities given in Lemma
\ref{lemaidentidadescoder} for  $n=3$. Such a proof can be
constructed in a similar way that the one for $n \geq 4$, taking
analog results for $\SUn(3)$-structures into account. Such
identities will be used in the following proposition, where we
compute some $\SUn(3)$-parts of $\nabla \psi_+$.
\begin{proposition} Let $M$ be a special almost Hermitian
$6$-manifold. Then
\begin{enumerate}
\item[{\rm (i)}]
 $\left( \nabla_{ \cdot } \psi_+ \right)_{1;a} = - w_1^+ e_i \otimes e_i \wedge \omega$,  $\left( d
\psi_- \right)_{1} =   2  w_1^+ \, \omega \wedge \omega$ and $
(d^* \psi_+)_{1} =  4 w_1^+ \omega$,
 where $w_1^+$ is
given by $ 12  w_1^+  = \ast ( d \psi_{-} \wedge \omega )=
\langle \ast d \psi_- , \omega \rangle$ or $\left( \nabla \omega
\right)_{1;+} = w_1^+ \, \psi_+$;
 \item[{\rm (ii)}]
 $\left( \nabla_{ \cdot } \psi_+ \right)_{1;b} = w_1^-  e_i \otimes Ie_i  \wedge \omega $,
$\left( d \psi_+ \right)_{1} = - 2  w_1^- \, \omega \wedge \omega$
and $(d^* \psi_-)_1 =  4 w_1^- \omega$, where $w_1^-$ is given by
$ - 12 w_1^- = \ast ( d \psi_{+} \wedge \omega ) = \langle \ast d
\psi_+ , \omega \rangle $ or $\left( \nabla \omega \right)_{1;-} =
w_1^- \, \psi_-$ ;
 \item[{\rm (iii)}] $ 4 \left( \nabla_{ \cdot } \psi_+ \right)_{1,2;a} = -  \langle \ast d
\psi_- , \omega \rangle \, e_i \otimes e_i \wedge \omega  +
\iota_{\omega}  \left( I_{(2)} - I_{(1)} \right) \ast d \psi_- $,
where $\iota_{\omega} \, : \, T^* M \otimes T^*M \to T^* M \otimes
T^*M \wedge \omega$ defined by $\iota_{\omega} ( a \otimes b ) = a
\otimes b \wedge \omega$;
 \item[{\rm (iv)}] $  - 4  \left( \nabla_{
\cdot } \psi_+ \right)_{1,2;b} = \langle \ast d \psi_+ , \omega
\rangle \, e_i \otimes I e_i \wedge \omega + \iota_{\omega}
 \left( 1 + I_{(1)} I_{(2)} \right)  \ast d \psi_+$
and $ \;- 2 \left( d \psi_+ \right)_{1,2} = - \langle \ast d
\psi_+ , \omega \rangle \, \omega \wedge \omega + \omega \wedge
\left( 1 + I_{(1)} I_{(2)} \right) \ast d \psi_+ $ ;
 \item[{\rm (v)}] $ 2 \left( \nabla_{ \cdot } \psi_+
\right)_{3} =  \Xi_+ \left( ( 1 - I_{(2)} I_{(3)} )  (d
\omega)_{3} \right)$, where $(d \omega)_{3} = (d \omega)_{3,4} -(d
\omega)_{4}$ with $4 (d \omega)_{3,4} = 3 d \omega + L(d \omega)$
and $\; 2 (d \omega)_{4} = - I d^* \omega \wedge \omega$.
\end{enumerate}
\end{proposition}
\Dt For part (i). If $(\nabla \omega)_{1;+} = w_1^+ \psi_+$,
$w_1^+ \in {\Bbb R}$, by Theorem \ref{igual3}, we obtain $ (\nabla
\psi_+)_{1;a} = \Xi_+ (\nabla \omega)_{1;+} =- w_1^+ e_i \otimes
e_i \wedge \omega$ and $(\nabla \psi_-)_{1;b} = \Xi_- (\nabla
\omega)_{1;+} = - w_1^+ e_i \otimes I e_i \wedge \omega$.
Therefore, $\left( d \psi_- \right)_{1} =  2 w_1^+  \omega \wedge
\omega$. On the other hand, since $ \omega \wedge \omega \wedge
\omega  = 6 \, Vol$, we have $ \ast \left( d \psi_- \wedge \omega
\right)  = \ast \left( \left( d \psi_- \right)_{1} \wedge \omega
\right) =  12 w_1^+=  \langle \ast d \psi_- , \omega \rangle$.

For part (ii). By an analog way, since $\left( \nabla \omega
\right)_{1;-} = w_1^- \, \psi_-$, $w_1^-  \in {\Bbb R}$,   we have
$ \left( \nabla_{ \cdot } \psi_+ \right)_{1;b} = \Xi_+ \left(
\nabla \omega \right)_{1;-}  = w_1^- e_i \otimes Ie_i \wedge
\omega$. Therefore, $ \left( d \psi_+ \right)_{1} = - 2 w_1^-
\omega \wedge \omega$. Hence we have $ \ast \left( d \psi_+ \wedge
\omega \right) = \ast \left( \left( d \psi_+ \right)_{1} \wedge
\omega \right) = - 12 w_1^-= \langle \ast d \psi_+ , \omega
\rangle$.

For part (iii). If $\nabla \psi_+ = x \otimes y \wedge \omega$, by
Theorem \ref{igual3},   we have
$$
4 (\nabla \psi_+)_{1,2;a} = x \otimes y \wedge \omega + y \otimes
x \wedge \omega + I x \otimes I y \wedge \omega + I y \otimes Ix
\wedge \omega.
$$
Therefore, $2 (d^* \psi_+)_{1,2} = - 2 \langle x, y \rangle \omega
+ Ix \wedge y - x \wedge Iy$. Since $ \langle d^* \psi_+ , \omega
\rangle = 12 w_1^+ = \, - 2 \langle x , y \rangle$ and
$
2 I_{(1)} (d^* \psi_+)_{1,2} = 2 \langle x,y \rangle \, \langle
\cdot , \cdot \rangle - ( x \otimes y + y \otimes x + Ix \otimes
Iy + I y \otimes Ix),
$
we have
$$
2 \iota_{\omega} I_{(1)} (d^* \psi_+)_{1,2} +   \langle d^* \psi_+
, \omega \rangle e_i \otimes e_i \wedge \omega = -4 ( \nabla
\psi_+ )_{1,2;a}.
$$
On the other hand, since $I (d^* \psi_+)_{1,2} = - (d^*
\psi_+)_{1,2}$ and $I(d^* \psi_+)_{4,5} = \linebreak - (d^*
\psi_+)_{4,5}$, it follows $2 (d^* \psi_+)_{1,2} = d^* \psi_+  + I
d^* \psi_+$. Thus,
$$ \iota_{\omega} (I_{(1)}- I_{(2)}) d^* \psi_+
+ \langle d^* \psi_+ , \omega \rangle e_i \otimes e_i \wedge
\omega =- 4  ( \nabla \psi_+ )_{1,2;a}.
$$
Finally, taking $ d^* \psi_+ = \ast d\psi_-$ into account, the
required identity in (iii) follows.

For part (iv). We proceed in a similar way as in the proof for
(iii), but now we consider
$$
4 (\nabla \psi_+)_{1,2;b} = x \otimes y \wedge \omega - y \otimes
x \wedge \omega + I x \otimes I y \wedge \omega - I y \otimes Ix
\wedge \omega
$$
and we compute $( d^*_{\omega} \psi_+)_{1,2}$. Thus we have $ 2 (
d^*_{\omega} \psi_+)_{1,2} = - 2 \omega(x,y) \omega + x \wedge y +
I x \wedge I y$. Since $(d^*_{\omega} \psi_+ )_1 = (d^* \psi_- )_1
=   4 w_1^- \omega = - \frac{2}{3} \omega(x,y) \omega$, we obtain
$$
2 \iota_{\omega} ( d^*_{\omega} \psi_+)_{1,2} + \langle
d^*_{\omega} \psi_+  , \omega \rangle  e_i \otimes I e_i \wedge
\omega = 4 (\nabla \psi_+)_{1,2;b}.
$$
Finally, taking $2 (d^*_{\omega} \psi_+)_{1,2} = 2 (d^*
\psi_-)_{1,2} = d^* \psi_-  + I d^* \psi_-$ and $d^* \psi_- = -
\ast d \psi_+$ into account, it follows the  first required
identity in (iv). By alternating both sides of such an identity,
the second required equation follows. Part (v) follows as in the
proof of Proposition \ref{calculoexolic}  for $(\nabla \psi_{+}
)_3$. \cqd

\begin{observation} {\rm
From the maps $\Xi_+$, $\Xi_-$ and identities (\ref{torsiones}),
it is not hard to prove
$$
\begin{array}{l}
\left( \nabla_{\cdot} \psi_- \right)_{1,2;a} = \left( \nabla_{I \cdot }
\psi_+ \right)_{1,2;b}, \qquad \left( \nabla_{\cdot} \psi_-
\right)_{1,2;b} = - \left( \nabla_{I \cdot } \psi_+ \right)_{1,2;a},  \\[2mm]
\left( \nabla_{\cdot} \psi_- \right)_{3,4} = - \left( \nabla_{I \cdot }
\psi_+ \right)_{3,4}.
\end{array}
$$
Thus, taking these identities into account, one can deduce the
respective ${\sl SU}(3)$-components of $\nabla \psi_-$ from those
of $\nabla \psi_+$. }
\end{observation}

 The following results are relative to  nearly K{\"a}hler six-manifolds.
\begin{theorem} \label{unomasunomenoscinco}
Let $M$ be a special almost Hermitian connected six-manifold of
type ${\cal W}_1^+ + {\cal W}_1^- + {\cal W}_5$ which is not of
type ${\cal W}_5$ such that $\nabla \omega = w_1^+ \psi_+ + w_1^-
\psi_-$, then
\begin{enumerate}
\item[{\rm (i)}] $\nabla \omega$ is nowhere zero, $\alpha =
(w_1^+)^2 + (w_1^-)^2$ is a positive constant and $d w_1^+ = -
w_1^- I \eta$, $d w_1^- = w_1^+ I \eta$;
 \item[{\rm (ii)}] the one-form $I\eta$ is closed and
given by $ 3 \alpha I \eta = w_1^+ d w_1^- - w_1^- d w_1^+$;
 \item[{\rm
(iii)}] $M$ is of type ${\cal W}_1^+ + {\cal W}_1^- $ if and only
if $w_1^+$ and $w_1^-$ are constant. \item[{\rm (iv)}] If $M$ is
of type ${\cal W}_1^+ + {\cal W}_5$, then $M$ is of  type ${\cal
W}_1^+$ or of type ${\cal W}_5$. \item[{\rm (v)}] If $M$ is  of
type ${\cal W}_1^- + {\cal W}_5$, then $M$ is of  type ${\cal
W}_1^-$ or of type ${\cal W}_5$.
\end{enumerate}
\end{theorem}
\Dt Since $M$ is of dimension six, it is straightforward to check
\begin{eqnarray}
\label{psipsimismo}&&(x \lrcorner \psi_+) \wedge \psi_+ =  (x
\lrcorner \psi_-) \wedge \psi_- = x \wedge \omega \wedge \omega =
- 2 \ast Ix, \\
\label{psipsidiferente} &&(x \lrcorner \psi_+) \wedge \psi_- =  -
(X \lrcorner \psi_-) \wedge \psi_+ = Ix \wedge \omega \wedge
\omega = 2 \ast x,
\end{eqnarray}
for all vector $x$.

 Since  $M$ is of type ${\cal
W}_1^+ + {\cal W}_1^- + {\cal W}_5$, we have
\begin{eqnarray}
d \omega & = & 3 w_1^+ \psi_+ + 3 w_1^- \psi_-, \label{nkdomega} \\
d \psi_+ & = & - 2 w_1^- \omega \wedge \omega  - 3 I \eta \wedge
\psi_-, \label{nkdpfimas}\\
d \psi_- & = &  2 w_1^+ \omega \wedge \omega  + 3 I \eta \wedge
\psi_+. \label{nkdpsimenos}
\end{eqnarray}
Now, differentiating equations
 (\ref{nkdpfimas}) and (\ref{nkdpsimenos}) and using equation (\ref{nkdomega}), we have
\begin{eqnarray}
\label{etauno}0 & = & 2 (d w_1^-  - 3 w_1^+ I \eta ) \wedge \omega
\wedge \omega + 3d I \eta \wedge
\psi_- , \\
\label{etados} 0 & = & 2 (d w_1^+ + 6 w_1^- I \eta)  \wedge \omega
\wedge \omega + 3 d I \eta \wedge \psi_+.
\end{eqnarray}
But $ d I \eta \in \Lambda^2 T^* M = {\Bbb R} \omega + {\frak
s}{\frak u}(3) + {\frak u}(3)^{\perp} $ and $ d I \eta_{{\frak
u}(3)^{\perp} } = x \lrcorner \psi_+$. Therefore,
$$
\begin{array}{c}
d I \eta \wedge \psi_+ = (x \lrcorner \psi_+ ) \wedge \psi_+ ,
\qquad \qquad  d I \eta \wedge \psi_- = (x \lrcorner \psi_+ )
\wedge \psi_- .
\end{array}
$$
Taking these identities into account and making use of equations
(\ref{psipsimismo}) and (\ref{psipsidiferente}), from equations
(\ref{etauno}) and (\ref{etados}) it follows
\begin{equation} \label{auxiluna}
\frac{3}{2} x = I d w_1^- + 3 w_1^+ \eta = - d w^+_1 - w_1^- I
\eta.
\end{equation}

On the other hand, differentiating equation (\ref{nkdomega}),
making use of equations (\ref{nkdpfimas}) and (\ref{nkdpsimenos}),
and taking $x \wedge \psi_+ = Ix \wedge \psi_-$ into account, we
obtain
$$
0 = (d w_1^+  + 3 w_1^- I \eta - I d w_1^-  - 3 w_1^+ \eta )
\wedge \psi_+.
$$
Therefore, taking equation (\ref{auxiluna}) into account,  we get
$I d w_1^- + 3 w_1^+ \eta = d w_1^+  + 3 w_1^- I \eta = 0$. Thus,
$d w_1^- = 3 w_1^+ I \eta$ and  $d w_1^+  = -  3 w_1^- I \eta$.
Moreover, $d \alpha = 2 (w_1^+ d w_1^+ + w_1^- d w^-_1) = 0$.
Since $M$ is connected, if $\alpha \neq 0$ in some point, then
$\alpha \neq 0$ everywhere. Now, it is immediate to check $3
\alpha I \eta = w^+_1 d w_1^- - w_1^- d w_1^+$  and $3 \alpha d I
\eta = 2 d w_1^+ \wedge d w_1^-=0$. Thus, parts (i) and (ii) of
Theorem are already proved.

Parts (iii), (iv) and (v) are immediate consequences of parts (i)
and (ii). \cqd

\begin{observation} {\rm
In \cite{nkstrGray},  Gray proved that if $M$ is a connected
nearly K{\"a}hler six-manifold (type ${\cal W}_1$) which is not
K{\"a}hler, then $M$ is an Einstein manifold such that ${\sl Ric} = 5
\alpha \langle \cdot , \cdot \rangle$, where ${\sl Ric}$ denotes
the Ricci curvature. In \cite{Cabrera:curvature},  showing an
alternative proof of such Gray's result,  we make use of Theorem
\ref{unomasunomenoscinco}.
 }
\end{observation}

\vspace{5mm}

\noindent {\bf 3.2 Four dimensions}

\noindent Now, let us pay lead our attention to manifolds with
${\sl SU}(2)$-structure.
\begin{theorem}
 \label{igual2} Let $M$ be a special almost Hermitian four-manifold
 with K{\"a}hler form $\omega$ and complex volume form $\Psi =
\psi_+ + i \psi_-$. Then
$$
\begin{array}{l}
\nabla \psi_+ \in T^* M \otimes \omega +  T^* M \otimes \psi_-,
\qquad \nabla \psi_- \in T^* M \otimes \omega +  T^* M \otimes
\psi_+,
\end{array}
$$
and $\Xi_{\pm} ({\cal W}_2) = \Xi_{\pm} ({\cal W}_4) = T^* M \otimes
\omega $. In this case, the space ${\cal W}= {\cal W}_2 + {\cal
W}_4$ of covariant derivatives of $\omega$ also admits the
relevant ${\sl SU}(2)$-decomposition \linebreak $ {\cal W} = T^*M
\otimes \psi_+ + T^*M \otimes \psi_- $, being ${\sl ker} \,
\Xi_{+} = T^*M \otimes \psi_-$ and  ${\sl ker} \, \Xi_{-} =  T^*M
\otimes \psi_+$.
\end{theorem}

If we consider the one-forms  $\xi_+$ and $\xi_-$ defined by
$\nabla \omega = \xi_+ \otimes \psi_+ + \xi_- \otimes \psi_-$,
i.e., $ \xi_+ = \langle \nabla_{\cdot} \omega , \psi_+ \rangle$
and $\xi_- = \langle \nabla_{\cdot} \omega , \psi_- \rangle$. The
two decompositions of $\xi$ are related as follows:
\begin{enumerate}
\item[(i)] $\xi \in {\cal W}_2$ if and only if $\xi_+ = I \xi_-$.
\item[(ii)] $\xi \in {\cal W}_4$ if and only if $\xi_+ = - I
\xi_-$.
\end{enumerate}
Moreover, we have the following consequences of last Theorem.
\begin{corollary} \label{dpsi+2}
For ${\sl SU}(2)$-structures, the exterior derivatives of
$\psi_+$, $\psi_{-}$  and $\omega$ are such that
$$
\begin{array}{l}
d \psi_+ = - \xi_+ \wedge \omega - 2 \eta \wedge \psi_+ =
 (\xi_+ \lrcorner \psi_- - 2 \eta) \wedge \psi_+,\\[2mm]
d \psi_- = - \xi_- \wedge \omega - 2 \eta \wedge \psi_- = -(
\xi_- \lrcorner \psi_+  + 2 \eta) \wedge \psi_- , \\[2mm]
d \omega = (\xi_+ - I \xi_-) \wedge \psi_+ = (\xi_+ \lrcorner
\psi_- - \xi_- \lrcorner \psi_+)  \wedge \omega.
\end{array}
$$
\end{corollary}
Hence the one-forms $\xi_+$, $\xi_-$ and $\eta$ satisfy
$$
\begin{array}{l}
  - \xi_+ \lrcorner \psi_- + 2 \eta = \ast \left( \ast d \psi_+ \wedge \psi_+
\right), \qquad   \xi_- \lrcorner \psi_+  + 2 \eta = \ast \left(
\ast d \psi_- \wedge \psi_-
\right), \\[2mm]
 \xi_- \lrcorner \psi_+ -  \xi_+ \lrcorner \psi_-  =   \ast \left( \ast d \omega \wedge \omega \right).
\end{array}
$$
Therefore, by Lemma \ref{estrellas}, we have
$$
\begin{array}{l}
4 \eta = \ast \left( \ast d \psi_+ \wedge
 \psi_+ \right) + \ast \left( \ast d \psi_- \wedge \psi_- \right)-
 \ast \left( \ast d \omega \wedge \omega \right),   \\[2mm]
2   \xi_- \lrcorner \psi_+   = \ast \left( \ast d \psi_- \wedge
\psi_- \right)   - \ast \left( \ast d \psi_+ \wedge
 \psi_+ \right) +\ast \left( \ast d \omega \wedge \omega \right),  \\[2mm]
2   \xi_+ \lrcorner \psi_-   = \ast \left( \ast d \psi_- \wedge
\psi_- \right)   - \ast \left( \ast d \psi_+ \wedge
 \psi_+ \right)-\ast \left( \ast d \omega \wedge \omega \right).
\end{array}
$$
Thus we can conclude that all the information about an
$\SUn(2)$-structure is contained in $d \omega$,  $d \psi_+$ {\it
and} $d \psi_-$ . Moreover, from these identities,  the equalities
for  $n=2$ contained in Theorem \ref{torsionw5} follow.

\vspace{2mm}

By Proposition \ref{cambio4}, for conformal changes of metric
given by $\langle \cdot ,\cdot \rangle_o = e^{2f} \langle \cdot,
\cdot \rangle$, we have $Id^* \omega_o = I d^* \omega - 4 df$ and
$\eta_o = \eta - 1/2 df$. The one-forms $\xi_+$ and $\xi_-$ are
modified in the way $ \xi_{+o} = \xi_+ - df \lrcorner \psi_-$,
$\xi_{-o} = \xi_- + df \lrcorner \psi_+$, where $\xi_{+o}$ and
$\xi_{-o}$   are the respective one-forms corresponding to the
metric $\langle \cdot , \cdot \rangle_o$. In fact, such identities
can be deduced taking the expression $ 2\nabla_o \omega_o = e^{2f}
\left\{ 2 \nabla \omega - e_i \otimes e_i \wedge I d f - e_i
\otimes Ie_i \wedge df \right\} $ into account,  where $\nabla_o$
is the Levi-Civita connection of $\langle \cdot , \cdot
\rangle_o$.

 \vspace{5mm}

 \noindent {\bf 3.1 Two dimensions}

\noindent Finally,  let us consider special almost Hermitian
two-manifolds. For these manifolds we have $\nabla \omega = 0$.
Therefore,
$$
\begin{array}{rcccl}
\nabla \psi_+ &=& - I \eta \otimes \psi_- & = & - \eta_+ \psi_-
\otimes \psi_- + \eta_- \psi_+
\otimes \psi_- \in {\Bbb R} + {\Bbb R},\\[2mm]
\nabla \psi_- &=& I \eta \otimes \psi_+ & = &  \eta_+ \psi_-
\otimes \psi_+ - \eta_- \psi_+ \otimes \psi_+ \in {\Bbb R} + {\Bbb
R},
\end{array}
$$
where $\eta = \eta_+ \psi_+ + \eta_- \psi_-$. Furthermore, $ d
\psi_+  =  - \eta_- \omega \in {\Bbb R} \omega$ and  $ d \psi_- =
\eta_+ \omega \in {\Bbb R} \omega$. Consequently,
 $ \eta_+ = -\ast d \psi_-$ and $\eta_- = \ast d \psi_+$.

With respect to  the curvature,  if $K$ denotes the sectional
curvature, it can be checked
$$
K(\psi_+ , \psi_-) = d I \eta (\psi_+ , \psi_-) = d \eta_+
(\psi_+) + d \eta_- (\psi_-) - \eta_+^2 - \eta_-^2.
$$

For conformal changes of metric given by $\langle \cdot , \cdot
\rangle_o = e ^{2f} \langle \cdot , \cdot \rangle$, the intrinsic
${\sl SU}(1)$-torsion is modified in the way $e^f \eta_{+o} =
\eta_+ - df(\psi_+)$ and $e^f \eta_{-o} =\eta_- - df(\psi_-)$,
i.e., $\;\; \eta_o = \eta - df $.

\begin{observation}
{\rm Let us consider an special almost Hermitian $2n$-manifold, $n
\geq 2$, which is K{\"a}hler (type ${\cal W}_5$). In such manifolds we
have
$$
d \psi_+ = - n \eta \wedge \psi_+ = - n I\eta \wedge \psi_-,
\qquad d \psi_- = - n \eta \wedge \psi_- =  n I\eta \wedge \psi_+.
$$
By differentiating these identities, it follows $ d \eta \wedge
\psi_+ =  d \eta \wedge \psi_- = 0$ and $d I\eta \wedge \psi_+ = d
I\eta \wedge \psi_- = 0. $ Therefore, $d \eta, d I \eta \in {\frak
s}{\frak u}(n) + {\Bbb R} \omega$.
 }
\end{observation}
\vspace{4mm}

\seccion{Almost hyperhermitian geometry}\indent
\setcounter{proposition}{0} A $4n$-dimensional manifold $M$  is
said to be {\it almost hyperhermitian}, if $M$ is equipped with a
Riemannian metric $\langle \cdot,\cdot \rangle$ and three almost
complex structures $I, J, K$ satisfying $I^2 = J^2= -1$ and $K= IJ
= - JI$, and $\langle A X, A Y \rangle = \langle X,Y \rangle $,
for all $X,Y \in T_x M$ and $A =I,J,K$. This is equivalent to
saying that $M$ has a reduction of its structure group to
$\SP(n)$. As it was pointed out in Section 2, each fibre $T_m M$
of the tangent bundle  can be consider as complex vector space,
denoted $T_m M_{\Bbb C}$,  by defining $i x = Ix$.

On $T_m M_{\Bbb C}$, there is an $\SP(n)$-invariant complex
symplectic form $\varpi_{I\Bbb C}= \omega_J + i \omega_K$ and a
quaternionic structure map defined by $y \to Jy$. Taking our
identification of $\overline{TM}_{\Bbb C}$ with $T^*M_{\Bbb C}$,
$x \to \langle  \cdot, x \rangle_{\Bbb C} = x_{\Bbb C}$, into
account (we recall $\langle \cdot , \cdot \rangle_{\Bbb C} =
\langle \cdot , \cdot \rangle + i \omega_I(\cdot, \cdot)$), it is
obtained $\varpi_{I\Bbb C} = Je_{i{\Bbb C}} \wedge e_{i{\Bbb C}}$,
where $e_{1}, \dots , e_{n}, Je_{1}, \cdots, Je_{n}$ is a unitary
basis for vectors. Therefore,
$$
\varpi_{I \Bbb C}^n = (-1)^{n(n+1)/2} n! \, e_{1{\Bbb C}} \wedge
\dots \wedge  e_{n{\Bbb C}} \wedge Je_{1{\Bbb C}} \wedge \dots
\wedge Je_{n{\Bbb C}}.
$$
Hence, we can fix $\Psi_I= \psi_{I+} + i \psi_{I-}$, defined by
$(-1)^{n(n+1)/2} n! \, \Psi_I = \varpi_{I\Bbb C}^n$, as complex
volume form.

By cyclically permuting the r{\^o}les of $I$, $J$ and $K$ in the
above considerations, we will obtain two more complex volume forms
$\Psi_J$ and $\Psi_K$. Thus, $M$ is really equipped with three
$\SUn(2n)$-structures, i.e., the almost complex structures $I$,
$J$ and $K$, the complex volume forms $\Psi_I$, $\Psi_J$, and
$\Psi_K$ and the common metric $\langle \cdot,\cdot\rangle$. We
could say that $M$ has a {\it special almost hyperhermitian}
structure. Furthermore, we also have
$$
(-1)^{n(n+1)/2} (n-1)! \, d \Psi_I = ( d \omega_J + i d \omega_K)
\wedge  (\omega_J + i \omega_K)^{n-1}.
$$
Hence,  we can compute $d \psi_{I+}$ and $d\psi_{I-}$ from $d
\omega_J$ and $d \omega_K$. Likewise, making use of considerations
contained in sections 2 and 3,   $\nabla \omega_I$ can be computed
from $d \omega_I$, $d \psi_{I+}$ and $d \psi_{I-}$. By a cyclic
argument, the same happens for $ \nabla \omega_J$ and $ \nabla
\omega_K$.
\begin{theorem} In an almost hyperhermitian manifold, the
covariant derivatives $\nabla \omega_I$, $\nabla \omega_J$ and
$\nabla \omega_K$ of the K{\"a}hler forms and the covariant derivative
$\nabla \Omega = 2 \sum_{A=I,J,K} \omega_A \wedge \nabla \omega_A$
are determined by the exterior derivatives $d \omega_I$,
$d\omega_J$ and $d\omega_K$.
\end{theorem}
In other words, $d \omega_I$, $d\omega_J$ and $d\omega_K$ contain
all the information about the intrinsic torsion of an ${\sl
Sp}(n)$-structure and the intrinsic torsion, determined by $\nabla
\Omega$ ( \cite{Swann:symplectiques,MartinCabrera}), of the
underlying ${\sl Sp}(n){\sl Sp}(1)$-structure. In relation with
last Theorem, we recall  Swann's result \cite{Swann:symplectiques}
that, for $4n \geq 12$, all the information about the covariant
derivative $\nabla \Omega$ is contained in the exterior derivative
${\sl d} \Omega = 2  \sum_{A=I,J,K} \omega_A \wedge d \omega_A$.
Furthermore, one of the consequences of previous Theorem is the
Hitchin's result \cite{Hitchin:Riemann-surface} that if the three
K{\"a}hler forms $\omega_I$, $\omega_J$ and $\omega_K$ of an almost
hyperhermitian manifold are all closed, then they are covariant
constant.  Almost hyperhermitian manifolds with covariant constant
K\"ahler forms  are called {\it hyperk\"ahler} manifolds. Such
manifolds are Ricci-flat. \vspace{2mm}

If the two  almost Hermitian structures determined by $I$ and $J$
are locally conformal K{\"a}hler (type ${\cal W}_4)$, then the one
determined by $K$ is also locally conformal K{\"a}hler
\cite{CabreraSwann}. Furthermore, in such a case, the three
structures have common Lee form. We recall that the Lee form is
defined by $\theta_A = - 1/(2n-1) A d* \omega_A$, $A=I,J,K$
\cite{GrayHervella}. Therefore, in such a situation we really have
a {\it locally conformal hyperk{\"a}hler} manifold. Let us compute
 the intrinsic torsion of the $\SUn(2n)_A$-structures, $A=I,J,K$.
 For $A=I$, we get
 $$
d \Psi_I =\frac{1}{ (-1)^{n(n+1)/2} (n-1)!} \theta \wedge
(\omega_J + i \omega_K)^n = n \theta \wedge \Psi_I,
 $$
where $\theta = \theta_I = \theta_J = \theta_K$. Therefore, $ d
\psi_{I+} = n \theta \wedge \psi_{I+} $ and, by Theorem
\ref{torsionw5}, we obtain that the ${\cal W}_5$-part of the
torsion is determined by
$$
\eta_I = \frac{1}{2n (2n-1)} I d^* \omega_I = - \frac{1}{2n}
\theta.
$$
Proceeding in a similar way for $J$ and $K$, we obtain $\eta_I =
\eta_J = \eta_K$. Furthermore, note that  the relevant one-form
$2n (2n-1) \eta_I - Id^* \omega_I$, given by Proposition
\ref{cambio4}, vanishes. In summary, we have the following result.
\begin{theorem} \label{example:aqh}
For  a locally conformal  hyperk{\"a}hler manifold of dimension $4n$
and a non null  Lee-form $\theta$, the three $\SUn(2n)$-structures
are of type ${\cal W}_4 + {\cal W}_5$. Moreover, the ${\cal
W}_5$-part of each one of such structures is determined by the
same one form $ \eta = - 1/2n \; \theta $.
\end{theorem}
As consequences of this Theorem, we have some results relative to
hyperk{\"a}hler manifolds.
\begin{corollary} \label{hiperka}
\begin{enumerate}
\item[{\rm (i)}] If the three $\SUn(2n)$-structures  of an almost
hyperhermitian $4n$-manifold are of type ${\cal W}_4$, then the
manifold is hyperk{\"a}hler.
 \item[{\rm (ii)}] For hyperk{\"a}hler manifolds, the intrinsic torsion of each
$\SUn(2n)$-structure vanishes.
\end{enumerate}
\end{corollary}

\begin{observation}{\rm Special almost Hermitian
manifolds with zero intrinsic torsion can be  called {\it ${\sl
SU}(n)$-K{\"a}hler manifolds}. The metric of such manifolds is
Ricci flat. Thus, Corollary \ref{hiperka} is an alternative proof
of the Ricci flatness of the hyperk{\"a}hler metrics.}
\end{observation}

{\small

\end{document}